\documentclass[12pt]{amsart}
\usepackage{epsf}

\makeatletter
\let\over\@@over
\let\Bbb\mathbb
\let\cal\mathcal
\makeatother

\newtheorem{lemma}{Lemma}
\newtheorem{proposition}{Proposition}
\newtheorem{theorem}{Theorem}
\newtheorem{remark}{Remark}
\newtheorem{claim}{Claim}

\def\Re{\mathop{\mathrm{Re}}}
\def\interior{\mathaccent "7017 }
\def\begfig {
\begin{figure}
\small
}
\def\endfig {
\normalsize
\end{figure}
}

\def\cases#1{\left\{\,\vcenter{\normalbaselines
    \ialign{$##\hfil$&\quad{##}\hfil\crcr#1\crcr}}\right.}

\begin{document}

\title{The Singly Periodic Genus-One Helicoid}

\author{David Hoffman}

\address{\hskip-\parindent
David Hoffman\\
Mathematical Sciences Research Institute\\
			1000 Centennial Drive  \\
			Berkeley CA 94720 	}

\email{david@msri.org}

\author{Hermann Karcher}

\address{\hskip-\parindent
Hermann Karcher			       \\
Mathematisches Institut Universit\"at Bonn		       \\
Wegerlerstrasse 8\\
5300 Bonn \\
Germany}

\email{unm416@ibm.rhrz.uni-bonn.de}

\author{Fusheng Wei}

\address{\hskip-\parindent
Department of Mathematics                \\
Virginia Tech\\                         
Blacksburg, VA  24061-0123}
\email{fwei@calvin.math.vt.edu}

\thanks{Hoffman was supported by research grant DE-FG03-95ER25250 of
the Applied Mathematical Science subprogram of the Office of Energy
Research, U.S. Department of Energy.  Hoffman and Wei were supported
by research grant DMS-95-96201 of the National Science Foundation,
Division of Mathematical Sciences.  Research at MSRI is supported in
part by NSF grant DMS-90-22140.}

\begin{abstract}
We prove the existence of a complete, embedded, singly periodic
minimal surface, whose quotient by vertical translations has genus one
and two ends.  The existence of this surface was announced in our
paper in {\it Bulletin of the AMS}, 29(1):77--84, 1993.  Its ends in
the quotient are asymptotic to one full turn of the helicoid, and,
like the helicoid, it contains a vertical line.  Modulo vertical
translations, it has two parallel horizontal lines crossing the
vertical axis.  The nontrivial symmetries of the surface, modulo
vertical translations, consist of: $180^\circ$ rotation about the
vertical line; $180^\circ$ rotation about the horizontal lines (the
same symmetry); and their composition.
\end{abstract}

\maketitle

\section*{Introduction}
\begfig
\hspace{2.0cm}
\epsfxsize=3.0in 
\epsffile{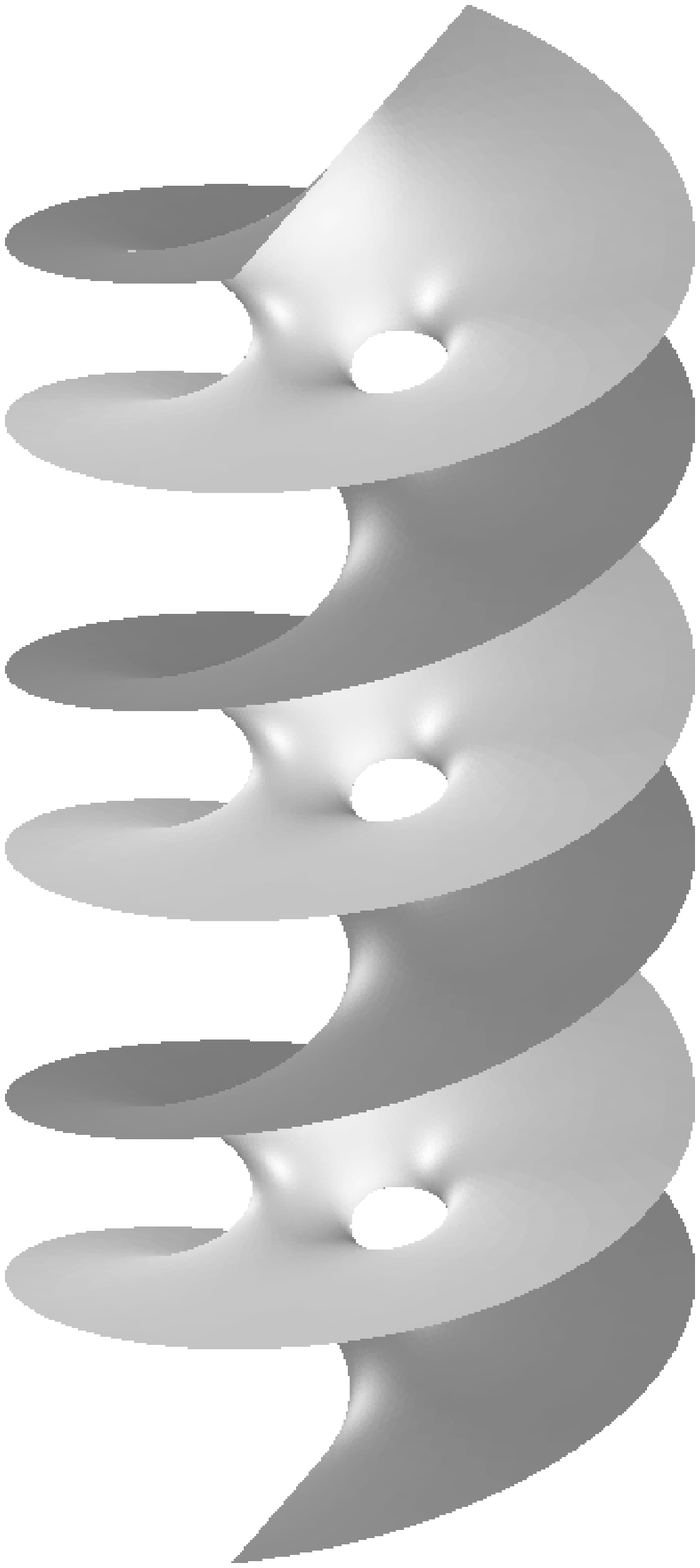} 
\hfill

\caption{The periodic genus-one helicoid.}
\normalsize
\endfig
In this paper, we prove the existence of a complete, embedded,
singly periodic minimal surface, whose quotient by vertical translations
has genus one and two ends.  The existence of this surface was announced
in \cite{howe2} 
and its significance discussed in \cite{howe3}.  
Its ends in the
quotient are asymptotic to one full turn of the helicoid, and, like the helicoid,
it contains a vertical line.  Modulo vertical translations, it has two parallel
horizontal lines crossing the vertical axis.  The nontrivial symmetries of
the surface, modulo vertical translations, consist of:
 $180^\circ$ rotation about
the vertical line; $180^\circ$ rotation about the horizontal lines (the
same symmetry); and their composition---a $180^\circ$ rotation about a line,
orthogonal to the lines on the surface, and passing through a common axis
point. This line meets the surface orthogonally and is referred to as a
{\em normal symmetry line}.

The description of the qualitative properties of the surface in the
paragraph above is sufficient to determine a two-parameter family of
Weierstrass data (1.7) that must contain the Weierstrass data
for any surface with these properties---if it exists. One parameter
controls the conformal type of the quotient, in this case a rhombic
torus. The other can be considered as controlling the placement of the
punctures corresponding to the ends. This is worked out in Section 1 and
presented in Theorem~1. 

The proof of existence of the singly periodic genus-one helicoid
consists of showing that the period problem ((1.8),(1.9)) is solvable.
This is done in Theorem 2 of Section 2. In  Theorem 3 of 
Section 3, we prove that the
surface is embedded by decomposing a fundamental domain into disjoint
graphs. As usual, the existence and embeddedness proofs
are independent. We do not use any special properties of the
parameters that kill the periods. In fact we show that {\it
any singly periodic (by translations) minimal surface that, in the quotient,
 is asymptotic
to the helicoid  (1.1) and contains a vertical axis and two horizontal
parallel lines  must be embedded}.

\bigskip

Other than the helicoid itself, this example was the first
 embedded minimal surface ever found that is asymptotic to the
 helicoid.
  It was one of the important steps in the discovery and construction
  of the non-periodic genus-one helicoid, whose existence is proved in
\cite{howe3}.  We hope that a complete understanding of this periodic
surface will be helpful in giving a non-computational proof, which is
not complete of this writing, of the embeddedness of the genus-one
helicoid.

We might have discovered this periodic surface earlier, had we been looking
 for it at the time.  In 1989, the first two authors realized that a
construction of Fischer and Koch \cite{fk1,fk2}
 could be modified to produce singly
periodic, embedded minimal surfaces with multiple helicoidal ends.
The Fischer-Koch triply periodic surface is formed of pieces congruent to the
solution to the disk-type Plateau Problem for the boundary in Figure 2.

\begfig
\centerline{\epsfxsize=1.4in 
\epsffile{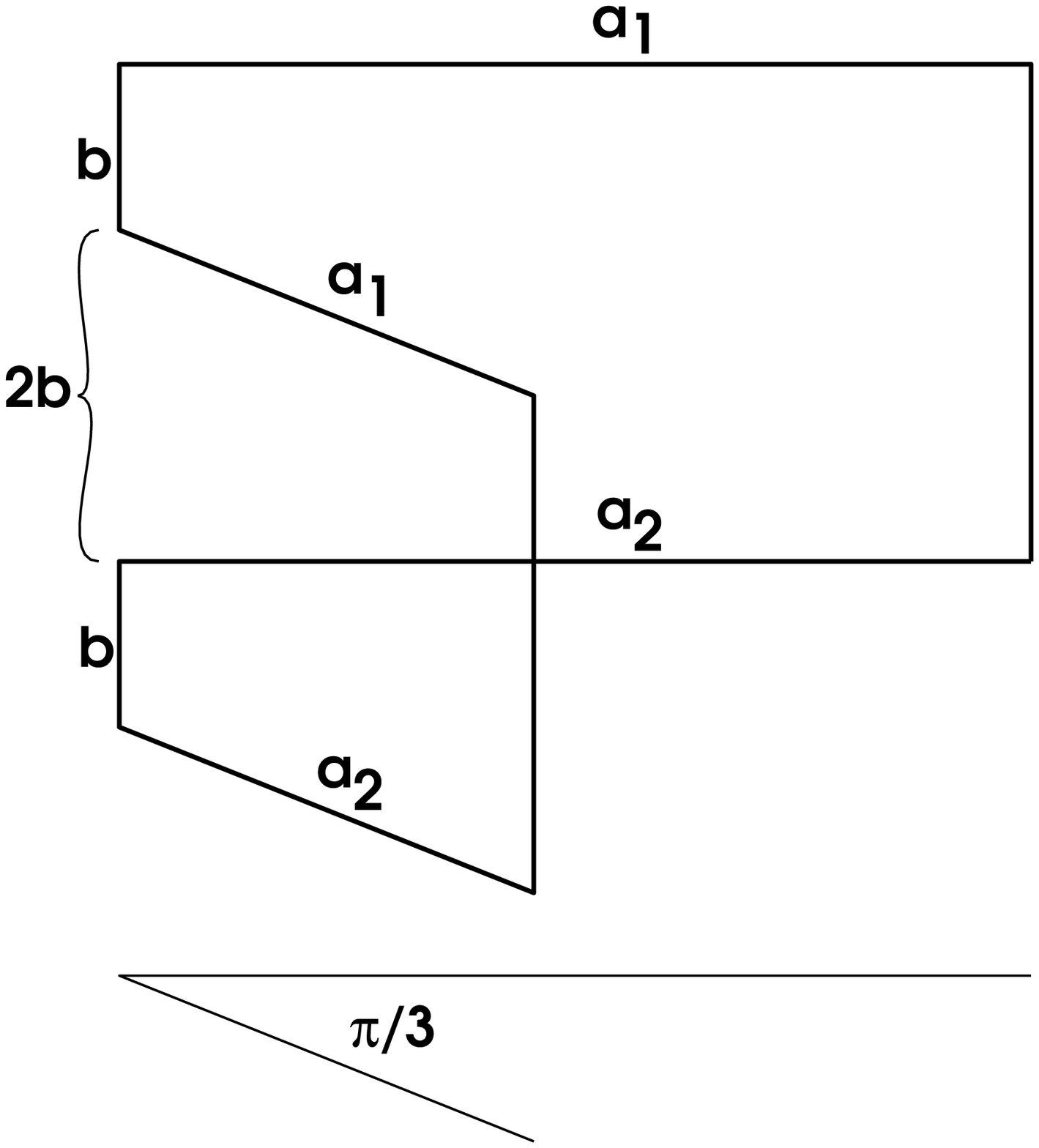}
\hspace{0.4cm}
\epsfxsize=2.8in 
\epsffile{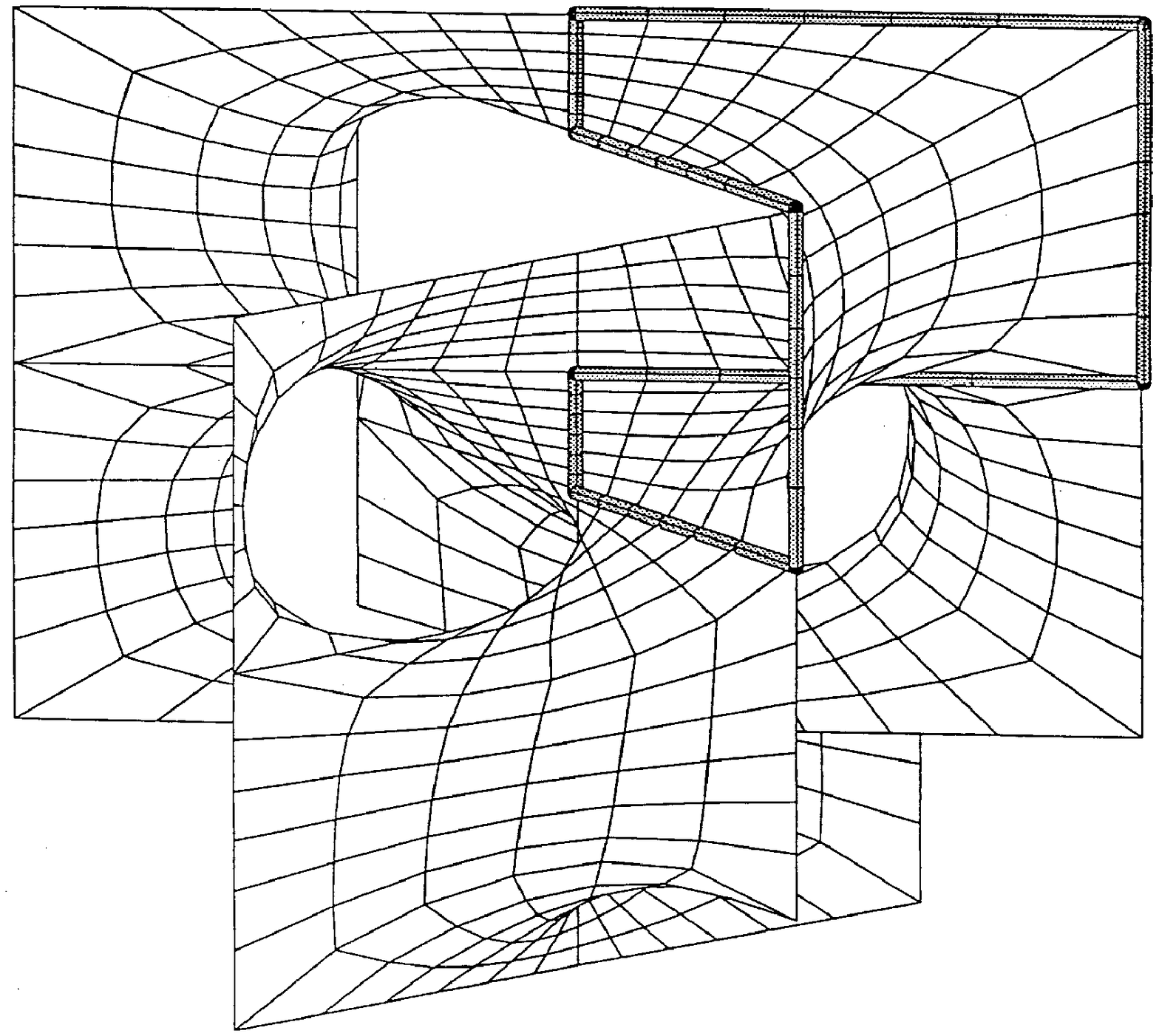}}

\vspace{0.2cm}

\caption{
{\em Above left}: The boundary of a fundamental building block of one
of the Fischer-Koch triply periodic surfaces.{\em Above right}: Six copies
of the fundamental building block,  which together form a
fundamental domain of the triply periodic surface. The surface is
computed using the discrete minimal surface ideas and  subsequent 
code of Pinkall
and Polthier \cite{pipo}. The computation
was carried out by Bernd Oberknapp. The original computation of the 
surface was done
by Ortwin Wohlrab; it was his pictures that suggested
to the first two authors the singly periodic surfaces described in the
text, one of which is illustrated below.   {\em Below}: A singly
periodic minimal surface, invariant under a vertical screw-motion and
asymptotic to three coaxial helicoids. It was suggested by the
Fischer-Koch surface, above right.}

\vspace{.5cm}

\hspace{4.2cm}
\epsfxsize=2.8in 
\epsffile{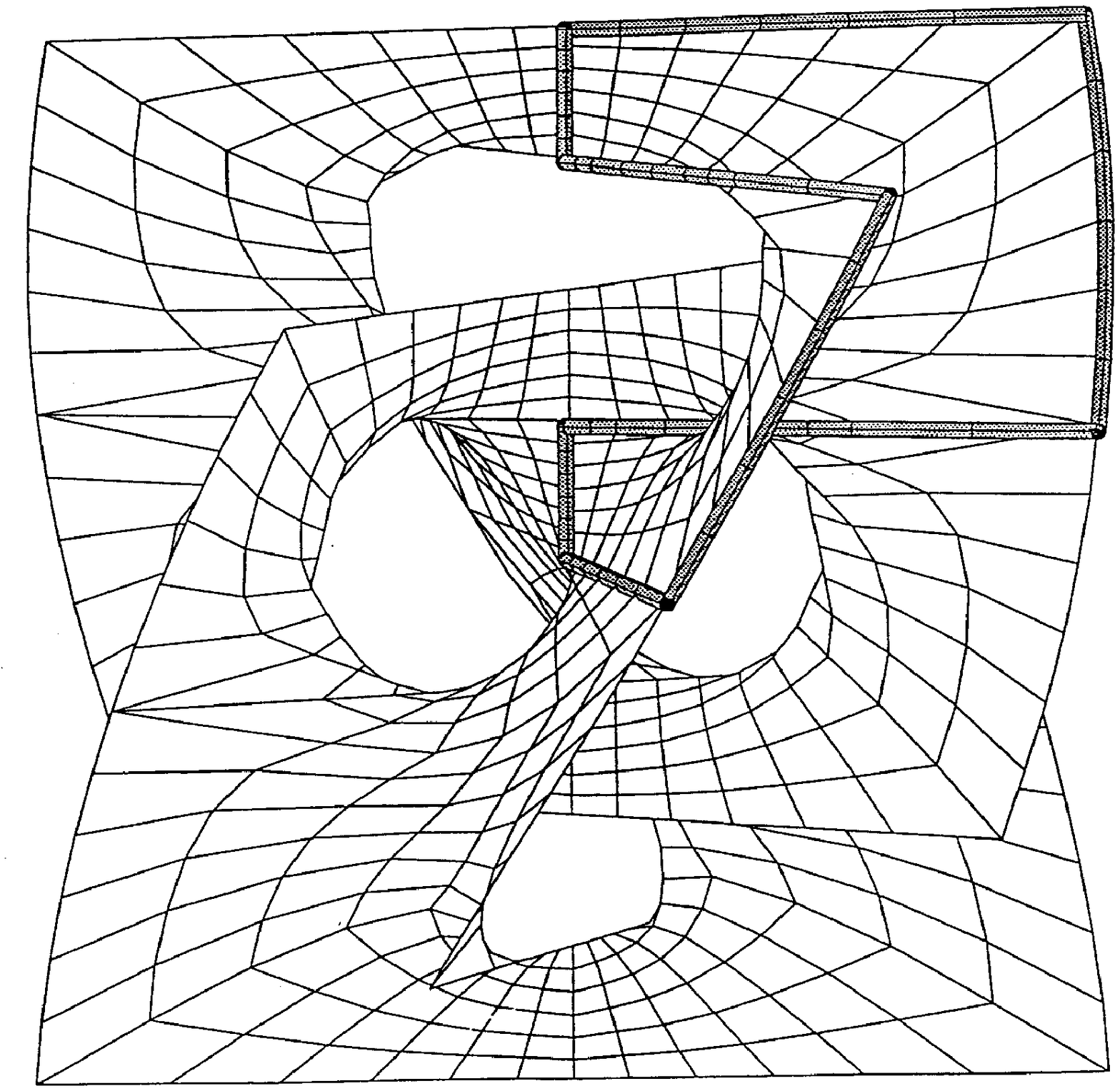}
\endfig

The surface extends, by $180^\circ$
rotation about its boundary line segments, to a triply periodic embedded
surface.  Our modification consisted of two simple steps.  First, we
realized that the length of the sides marked $a_i$ could be increased
without limit, producing an embedded minimal graph over a strip.  This
extends to an embedded, singly periodic surface with six flat ends of
Scherk-type.  Second, we observed that the fundamental piece could be
modified by rotating the horizontal sides $a_1$ and $a_2$ by a fixed
angle, say $\theta$.  Each $\theta$ produces a fundamental embedded piece
that extends by $180^\circ$ rotation to an embedded minimal surface,
asymptotic to three coaxial helicoids, and invariant under a vertical
screw motion of the form
$$p \to e^{2i\theta} p+ (0, 0, 8b).$$

Its Weierstrass representation is suggestive of the surface
 to which we now turn
our attention.

\section{Determination of the Weierstrass Representation}

\renewcommand{\theequation}{1.\arabic{equation}}

 The
helicoid can be described by the  data
\begin{equation} 
g = z, ~dh = {{i\,dz}\over z}
\end{equation}
on  ${\cal S} = {\Bbb C} - \{0 \}$ in the Weierstrass
representation 
\begin{equation}
 X (p) = X (p_0) + \Re \int_{p_0}^p \Phi, 
\quad \Phi = \left ( {1\over 2} (g^{-1} -g), ~{i\over 2} ( g^{-1} +g),
1\right) \,dh. 
\end{equation}
 The integration 
produces a conformal minimal immersion, but has a period equal to 
$(0,0, \pm 2\pi)$ on any closed curve $\gamma$ homotopic to
 $\vert z\vert =1$. Thus the immersion in (1.2) with data (1.1) is multivalued, 
and
its image, the helicoid, is invariant under a vertical translation
produced by the period of (1.2) on the cycle $\vert z\vert =1$.

The minimal surface we wish to construct is singly periodic, has a
helicoidal-type end and, modulo translations, will have genus equal to
one.  From general results \cite{mr2} and our assumptions about the
geometry and topology of the surface in the Introduction, we know that
(if it exists) it will have a Weierstrass representation on a
twice-punctured torus and that $g$ and $dh$ will extend meromorphically
to the compact torus.\footnote{The helicoid data in (1.1) is defined
on a twice-punctured sphere and it clearly extends meromorphically.}
Our goal is to determine the Weierstrass representation from this
information and from the symmetry we assume the surface to have.

\subsection{The lines on the surface}

We assume that, like the helicoid, our surface contains a single vertical
line and---again like the helicoid---the symmetry of $180^\circ$ rotation about
this line fixes this line and no other points.  This implies that the
torus is a rhombic torus, which we denote by $T^2$.

We assume in addition that the surface contains, modulo translations,
two parallel {\it horizontal} lines, each meeting the vertical line in
a single point.  Rotation of $180^\circ$ about either of these lines fixes
both lines in the quotient.  Because these lines diverge, 
such a rotation leaves the punctures corresponding
to the ends unmoved .  Each line
diverges to {\em both} ends, forcing the fixed-point set on the torus 
of this $180^\circ$ rotation to be a single closed symmetry curve that
must contain both end-punctures.  This curve must cross the curve
corresponding to the vertical axis in two points.  The two punctures
separate the symmetry line into two components corresponding to the two
horizontal lines.  These two points must be symmetrically placed with
respect to $180^\circ$ rotation about the vertical line.

We choose a fundamental domain for our torus so that the lines on the
surface correspond to the diagonals of a rhombus. Because we will use
our development of elliptic functions in \cite{howe3}, these diagonals
are placed to make a $\pm 45^\circ$ angle with the real axis in ${\Bbb C}$. 
 See Figure~4.  We think of them as the {\em horizontal} and the
{\em vertical} diagonal, according to the lines on the surface in
space, which are their images. To avoid possible confusion with lines
in the $\Bbb C$ and in $\Bbb R^3$, we will refer to them as the
{v-diagonal} and the {h-diagonal}. Without loss of generality, we may
assume that the h-diagonal is the one making a  $45^\circ$ angle with
the real axis.

Wherever possible, we will deal only with elliptic functions, avoiding
any notational reference to the complex variable in $\Bbb C$,  which is
not well defined on the torus. However we will have need of $du$---the
differential of a complex variable $u$---which projects to the torus and
is free of poles and zeros. 

 Rotation by $180^\circ$ about a point on the surface
where the straight lines meet is a rotation about the line determined
by the normal to the surface at that point, we refer to it as a {\em
normal symmetry}.  On the quotient torus it is an
orientation-preserving involution: $180^\circ$ rotation about the
center of the torus.  We denote this rotation by ${r\!_{_P}}$. It fixes
four points that are on a half period lattice. We may choose as a
fundamental domain for the torus a rhombus whose center is at one of
these points and whose vertex is another such point.  These two points
are the intersection of the diagonals of the rhombus.  The other two
fixed points of the normal symmetry both lie at the same height as one
or the other of the fixed points that lie on the vertical axis.  {\em
We assume, without loss of generality, that they lie at the same height
as the fixed point corresponding to the center of the rhombus.}

\begfig
\hspace{2.2cm}
\epsfxsize=3.0in 
\epsffile{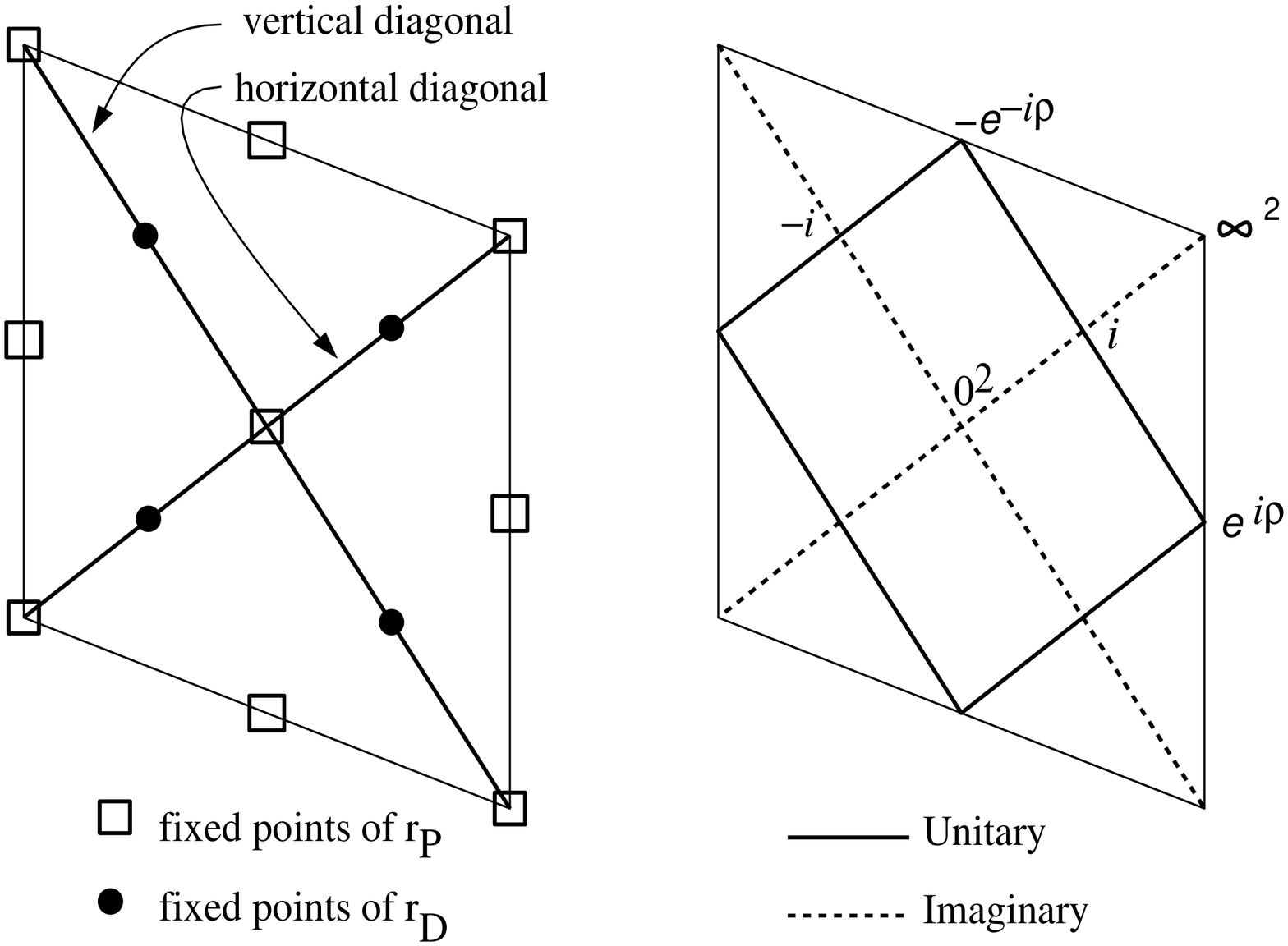} 
\hfill 
\vspace{0.2cm}

\caption{ {\em Left}: Fixed points of $r_D$ and $r_P$ The diagonals are
labeled vertical and horizontal here because they are the preimages of
the vertical and horizontal lines on the surface. In the text, they are
referred to as the v- and h-diagonals. {\em Right}: Values of the
elliptic function $z$. }

\endfig
\subsection{Elliptic functions associated to our construction}

Euler's formula implies that $T^2/{r\!_{_P}}$ is the sphere.  The
quotient map $T^2 \to T^2/{r\!_{_P}}$ is meromorphic, has degree two,
and is branched at the half-period points.  We want to use this map to
describe the torus analytically and for that it is necessary to turn
this map into a function by identifying three points of the sphere
$T^2/{r\!_{_P}}$ with points in the complex plane.  The motivation for
our choice is to make the symmetries of the torus induce simple
symmetries of $\Bbb C \cup \{ \infty \}$.  This is explained in detail
in \cite{howe3}.

Denote by $O$ (resp. $O^\prime$) the center point (resp. the vertex) of
$T^2$.  The choice we make is to have $0\in \Bbb C$ equal the
projection of $O$, and $\infty$ equal the projection of $O^\prime$
(since these points are fixed by ${r\!_{_P}}$, these values are branch
values) and $+i$ equal the projection of
 the midpoint of the h-diagonal.  This determines a degree-two
elliptic function, which we denote by $z$.  See Figure~3.

 Let ${r\!_{_D}}$ denote $180^\circ$ rotation
around a midpoint of a diagonal between $O$ and $O^\prime$. It fixes all
four such midpoints. Observe that ${r\!_{_D}}$ interchanges $O$ and $O^\prime$
in $T^2$, where $z = 0$, $\infty$ respectively, and fixes a point where $z =i$.
This implies that $z\circ {r\!_{_D}} = -1/z $.  Thus $z = \pm i$ at the
fixed points of ${r\!_{_D}}$.  Since $z \circ {r\!_{_P}} =z$, by definition, 
and the 
degree of $z$ is two, the values of $z = \pm i$ at the fixed points of
${r\!_{_D}}$ must be as in Figure~3.

Let $\mu$ be reflection in a diagonal of $T^2$.  Then $\overline{z \circ
\mu} = -z$.  To see this, first observe that both functions in question
are meromorphic on $T^2$. Their quotient has no zero or poles, and takes
on the value 1 at a fixed point of ${r\!_{_D}}$.  This implies that $z$ has
imaginary values on the diagonals.  Similarly, let $\nu$ denote reflection
in the line parallel to a diagonal and passing through a fixed point of
${r\!_{_D}}$.  Then ${\overline z\circ \nu} = 1/z$.  Thus $\vert z \vert
=1$ on this line. 

\begfig
\hspace{2.8cm}
\epsfxsize=2.6in 
\epsffile{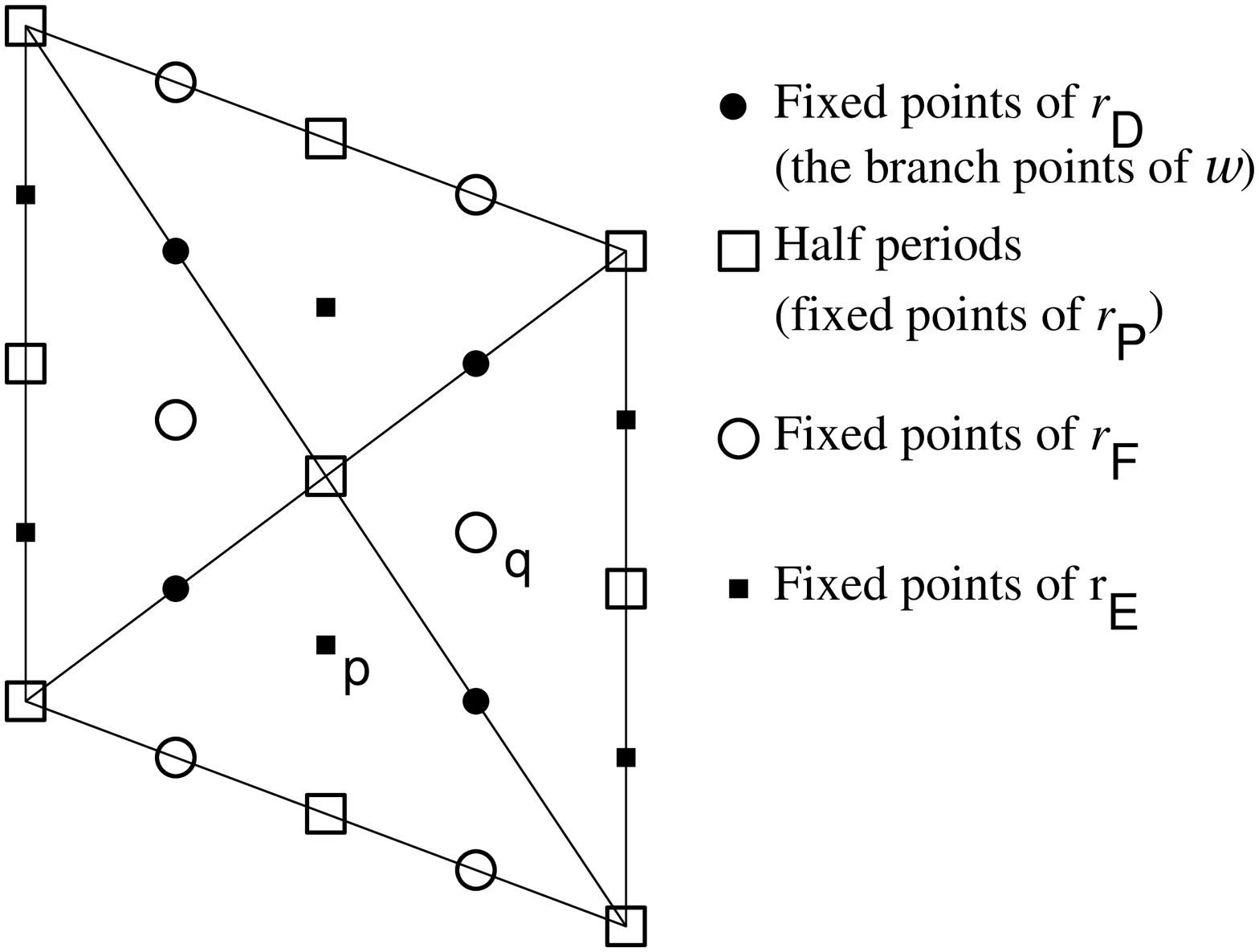} 
\hfill
\vspace{0.3cm}

\caption{
The fixed point sets of the four involutions of the torus. }

\endfig

In particular, the values of $z$ at the two half-period points other
than $O$, $O^\prime \in T^2$ are of the form $e^{i\rho}$ and 
$-e^{-i\rho}$, because
${r\!_{_D}}$ interchanges these points and $z \circ {r\!_{_D}} = -1/z$.

In fact, $z$
is a geometric normalization of the Weierstrass $\wp$-function.  It
satisfies the differential equation

\begin{equation}
({{z'}\over{z}})^2 = -{{2}\over{\cos \rho}}({z-{{1}\over{z}}-2i\sin \rho})
\end{equation}

 The value of $\rho$ characterizes the rhombic torus as
can be seen explicitly in (1.3) above and Lemma 1 iii) below. 

 Before returning to the determination of the Weierstrass data, we will
construct another elliptic function intimately related to the minimal
surface. This time we will use the involution ${r\!_{_D}}$, which fixes
four points on the diagonals.
 The projection $T^2\to T^2/{r\!_{_D}}$ is a degree-two meromorphic map
that we will turn into an elliptic function by specifying its values at
three points of $S^2 = T^2/{r\!_{_D}}$.  At the projection of $O \in T$
we specify the value $0$.  At the projection of one of the half periods
(not $O^\prime$) we specify the value $\infty$.  Because $w \circ
{r\!_{_D}} = w$ by definition, this determines $w$ at the half-period
points. (See Figure~5.)  We now determine $w$ by specifying its value to
be $+1$ at the point $q$ indicated in Figure~4. 

 The following Lemma
contains information we will need about $w$ in our discussion of the
surface. In the proof of the lemma
and in discussions that follow, we will need to distinguish between
reflection in the h-diagonal and  the v-diagonal. Let
$\mu$ denote reflection in the h-diagonal and  write
$\mu_{vert}$ for reflection in the v-diagonal.

\begfig
\hspace{2.5cm}
\epsfxsize=3.0in 
\epsffile{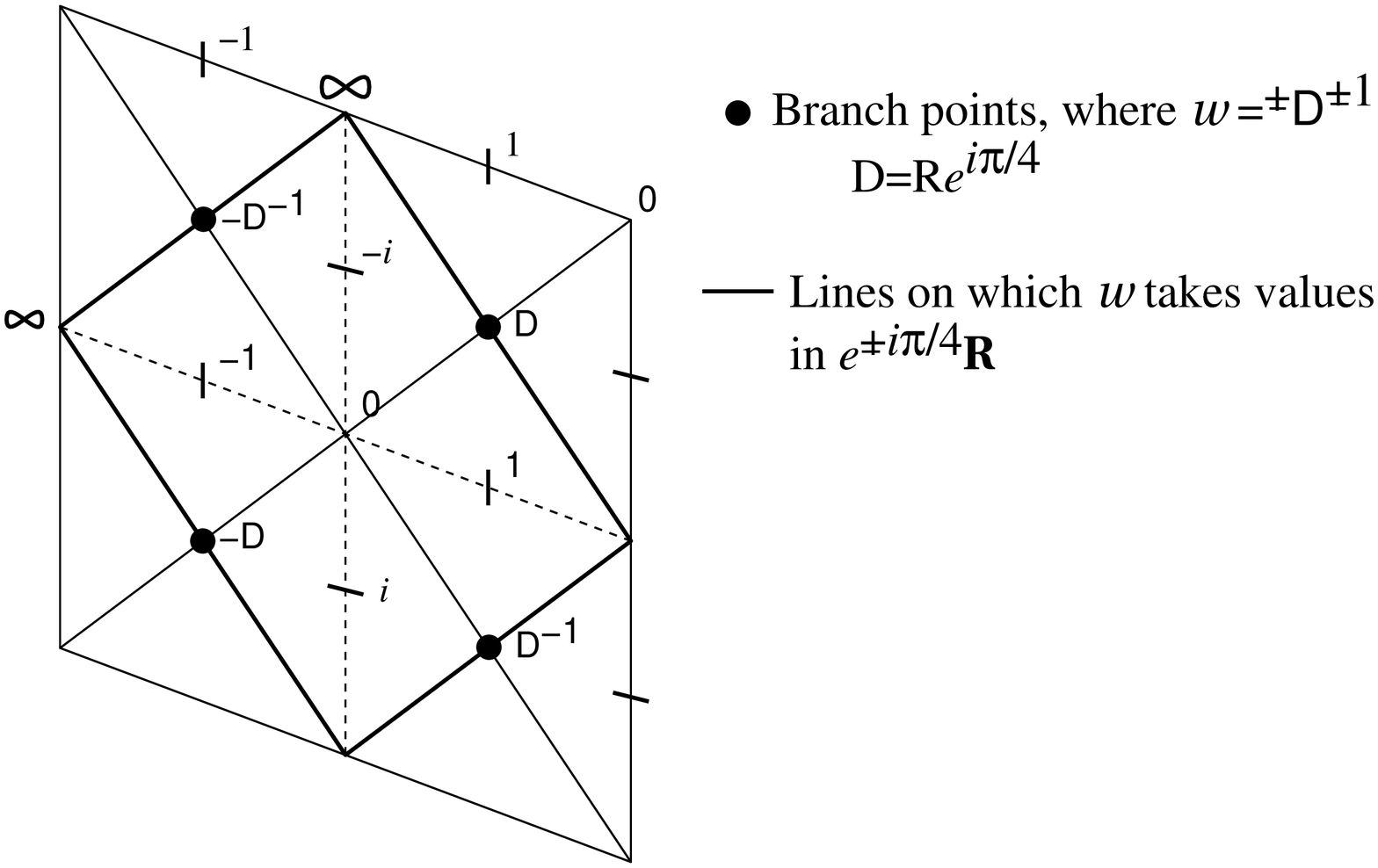} 
\hfill

\caption{The values of the elliptic function $w$.}

\endfig

\begin{lemma}
 The degree-two elliptic function $w$ on a rhombic torus has values as
indicated in Figure~5:

\begin{enumerate}
\item[(i)] Along the diagonals and along the lines through the branch points, $w$
takes on values in $e^{\pm i\pi/4} {\Bbb R}$.  In particular, its branch
values are of the form $\pm R e^{\pm i \pi/4}$, the constant $R > 0$ being
determined by the choice of rhombus;

\item[(ii)] For a fixed rhombic torus, the branch values of $w$ are related to
the branch value $e^{i\rho}$ of $z$ by $R^2 = \cot (\pi/4-\rho/2)$;

\item[(iii)] The rhombic torus is determined by the algebraic relation
$$(2\cos \rho) w^{-2} = - (z - z^{-1} - 2i \sin \rho).$$
\end{enumerate}
\end{lemma}

{\bf Proof of Lemma 1.} 
We will prove (i) here.  Statements (ii) and (iii) are proved in
\cite{howe3},
where there is a fuller discussion of Jacobian elliptic functions.

By definition, $w \circ {r\!_{_D}} =w$, so the placement of the points where
$w =1$ is determined.  Let $r_F$ denote $180^\circ$-rotation around $q$, an
involution whose other fixed points are marked by $\circ$ in Figure~4.  
The meromorphic functions $w \circ r_F$ and $w^{-1}$ agree at $q$ and have the
same zeros and poles.  Hence $w \circ r_F = w^{-1}$.  We already knew $w
=1$ at two of the fixed points of $r_F$; hence $w = -1$ at the other two.
In the same manner $w \circ {r\!_{_P}}$ and $w$ have the same zeros and poles and,
noting  the values of $w$ at the fixed points of $r_F$, one sees that
$w \circ {r\!_{_P}} = -w$.  Similarly, $r_E$ is $180^\circ$ rotation about $p$
(see Figure~4), and $w \circ r_E = -w^{-1}$.  This implies that $w = \pm i$
at the fixed points of $r_E$. The signs are determined by the
fact that $w$ is orientation preserving.

 The meromorphic
 function $\overline {w\circ \mu_{vert}}$  has the same zeros and
poles as $w$.
 Evaluation at $q$ gives the relation $\overline {w\circ \mu_{vert}}
= -iw$.  This means that $w$ takes on values in $e^{-i\pi/4} {\Bbb R}$ on
the v-diagonal.  If $ \mu$ is reflection in the h-diagonal
$\overline{w\circ\mu} = i w$ and therefore $w$ takes on values in $e^{i
\pi/4} {\Bbb R}$ on this line.  Similarly, reflection $\nu $ (resp.{$\hat\nu$})
in the lines through the other two fixed points of ${r\!_{_P}}$ satisfies 
$\overline {w\circ \nu} = iw$ (resp. $\overline {w \circ \hat \nu} = -iw$), so $w$ takes values in
$e^{i\pi/4} {\Bbb R}$ (resp. $e^{-i\pi/4} {\Bbb R})$ on these lines.

\begfig
\hspace{4.1cm}
\epsfxsize=1.2in 
\epsffile{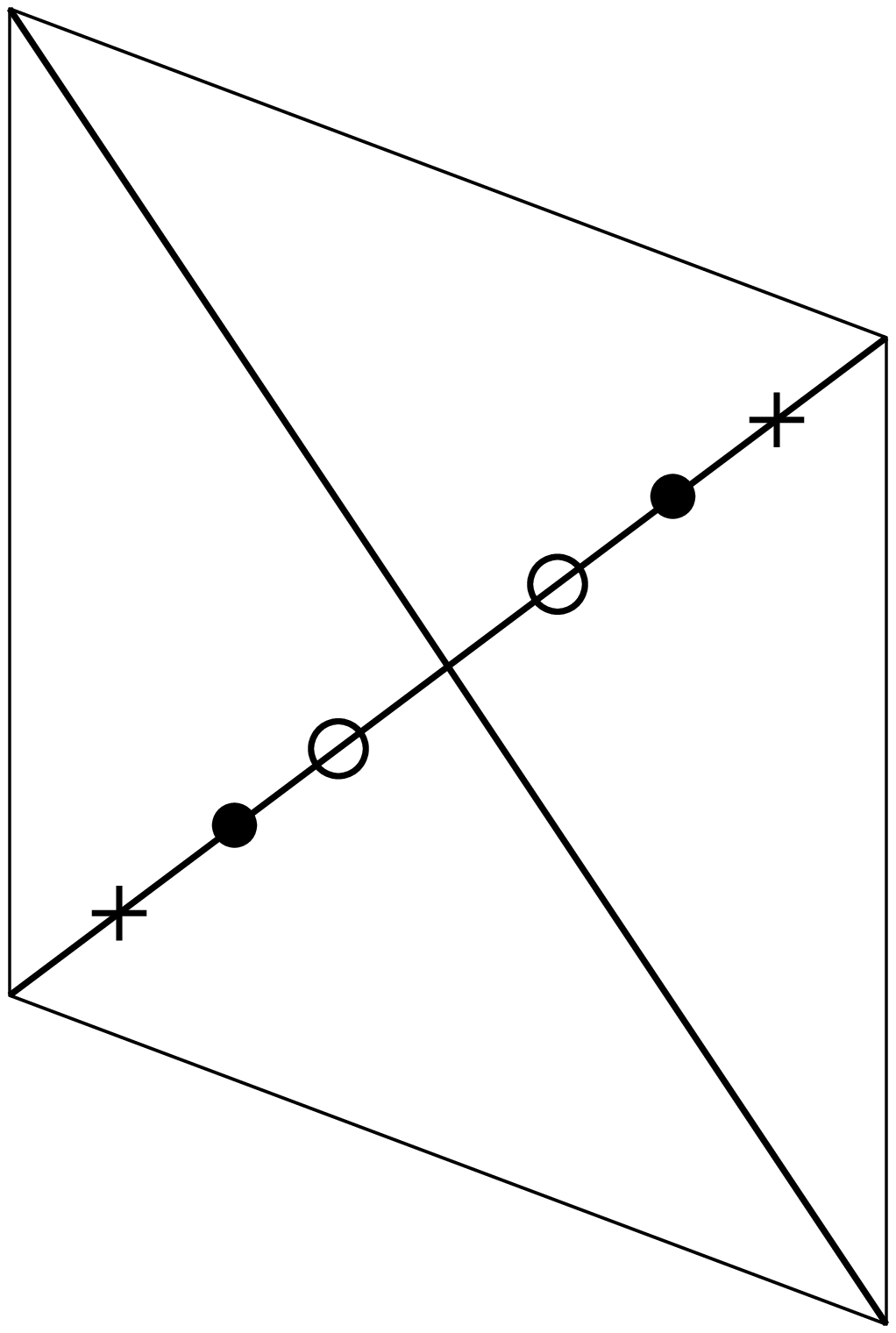} 
\hfill

\def\thecaption{ A $\bullet$ marks the branch points of $g$ on the
h-diagonal, and a $\circ$ marks the points on the surface where the
Gauss map is vertical, i.e. $g$ is zero or infinity. An x marks the
punctures. For this illustration, we are assuming that $\lambda <1$;
for $\lambda >1$, the punctures (x) and vertical points ($\circ$) are
 interchanged.} \caption{\thecaption}

\endfig


\subsection{The Gauss map in terms of $w$}

From \cite{mr2},
we know that the total curvature  of our desired surface is $-2\pi (\chi(M)- W(M))
$, where $\chi(M)$ is the Euler characteristic of the quotient surface $M$, and
$W(M)$ is the total winding number of the ends of $M$ . In our case, $M$
is a torus $T^2$ punctured twice, so $\chi(M)= -2$. Both ends are 
asymptotic to a single full turn of the helicoid, so each end contributes
1 to $W(M)$. Hence  the total curvature is $-8\pi$ and the degree of $g$ is 2.

Because we want a helicoidal end on the surface, as well as a vertical line  that
corresponds to the axis of the helicoid, we see
from (1.1) that $g =0$ at one puncture and $g = \infty$ at the other.
There must be one other zero and one other pole.  Because rotation about a
vertical (or horizontal) line in ${\Bbb R}^3$ will leave zeros and poles
of $g$ unchanged (or interchange them), the other zero and pole of $g$
must lie on the symmetry lines; otherwise there would be too many
 zeros and poles.  On the
v-diagonal, $g$ must be unitary, so the zeros and poles are on the
h-diagonal.

As in the case for any degree-two elliptic function, the branch points
are symmetrically placed with respect to the zeros and poles.(See
\cite{howe3}.)  Because the end-punctures are also on the h-diagonal,
two branch points are on this diagonal, too.  The four branch points of
$g$ form a set that is invariant under symmetries of the surface.  In
particular, the branch-point set is invariant under ${r\!_{_P}}$
(induced by the normal symmetry) and ${r\!_{_D}}$.  This implies that
the branch points of $g$ are the quarter-points of the diagonals; these
are exactly the branch points of $w$.  Thus $g$ and $w$ differ by a
M\"obius transformation.

 At a
puncture,  $w = \pm r e^{i\pi/4}$ for some $r > 0$, a consequence of
Lemma 1(ii) in Section 1.2.  Because $g$ is unitary at $O \in T^2$, it
follows that 
$$g = e^{i\theta} {{w-re^{i\pi/4}}\over {w +
re^{i\pi/4}}}$$
 and, after a rotation about a vertical axis, if
necessary,  we may assume that $\theta = 0$, i.e.
  \begin{equation}
 g =
{{w-re^{i\pi/4}}\over {w + re^{i\pi/4}}}. 
\end{equation}

\subsection{The  complex height differential in terms of $z$}

We can easily determine the differential $dh$, which is holomorphic on
the punctured surface. From (1.2)  (or see \cite{hk2}) the metric on
the surface is given by
 $ds = (\vert g \vert + \vert g \vert^{-1}) \vert dh\vert$.  Because we
require $ds$ to be everywhere nonzero (for regularity), $dh$ has simple
zeros at the two points where the Gauss map is vertical on the surface
.  Thus it has only two poles, which must be located at the punctures.
In particular, $dh \circ {r\!_{_P}} = dh$.


In Section 1.1 we set up the parameter domain so that
the  h-diagonal would be mapped into a horizontal line
in $R^3$. This forces $dh$ to be imaginary on
the h-diagonal; i.e.  $x_3 = \Re \int dh$ must be constant
on this line.

Just as $w$ was well adapted to $g$, $z$ is a good match for $dh$.
Recall that $z$ is imaginary on both diagonals and, by definition, 
 $z \circ {r\!_{_P}} = z$.  We define
$\lambda$ by the requirement that $z = i/\lambda$ at the punctures.  Since
$z \circ {r\!_{_D}} = -1/z$, $z = i \lambda$ at the vertical points of $g$ on the
h-diagonal.
Note that, by its definition, $\lambda >0$. Also, $\lambda \neq 1 $ 
because this $g$ is not branched at a helicoidal end.  The
function $${{z-i\lambda}\over {z-i\lambda^{-1}}}$$ has the same poles and
zeros as $dh$ and is real on the diagonals.

We can express $dh$ in terms of the standard pole-and-zero-free
holomorphic form $du$, described in Section~1.1.  Under our assumption
that the h-diagonal makes a $45^\circ$ angle with the real
axis, $e^{i\pi/4} du$ is imaginary on this diagonal. Hence, we may
assume without loss of generality that \begin{equation} dh = e^{i\pi/4}
{{z-i\lambda}\over {z-i\lambda^{-1}}} du, \end{equation} because we are
free to scale the surface by a real constant.

\begin{remark}.  Using Lemma 1 (iii) we can relate $\lambda$ to r: 
 
\begin{equation}
 \lambda + \lambda^{-1}= 2(\sin \rho + (\cos \rho) /r^2).
\end{equation}
\end{remark}

\subsection{The Ansatz for the periodic genus-one helicoid}

\begfig
\hspace{1.0cm}
\epsfxsize=4.2in 
\epsffile{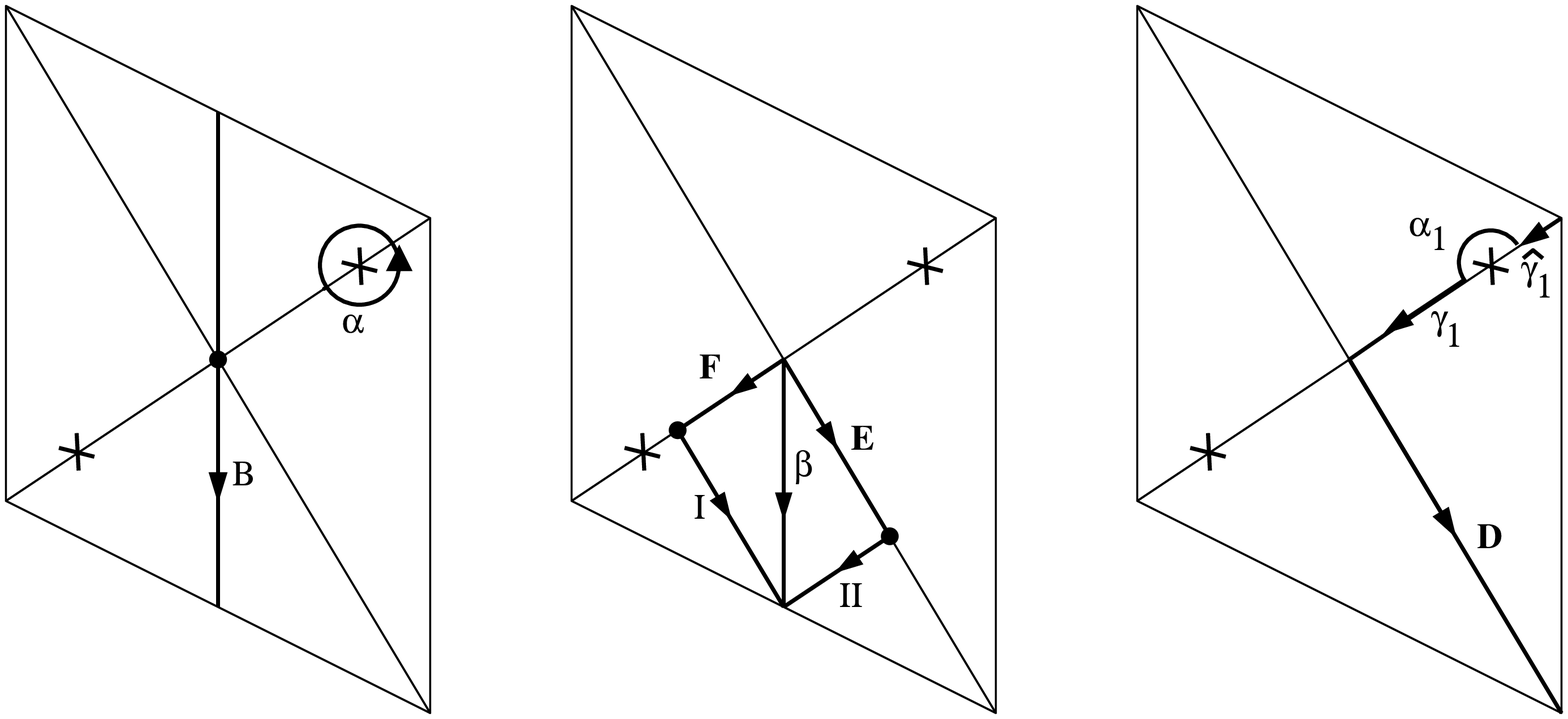} 
\hfill 
\vspace{0.2cm}

\def\thecaption{
{\em Left}: Two cycles that, together with their reflections, generate a 
homology basis of the punctured torus. {\em Center}: Paths used to derive the period
conditions (1.8) and (1.9).
{\em Right}: A cycle homologous to $B$.}
\caption{\thecaption}
\endfig

Our geometric description of the surface has led to an explicit
two-parameter family of Weierstrass data: For $\rho \in (-\pi/2, \pi/2)$ and
$\lambda \in (0,1)\cup (1, \infty)$,
 \begin{eqnarray}
 g =
{{w-re^{i\pi/4}}\over {w + re^{i\pi/4}}},
\mbox{\,\,\,\,\,\,} dh = e^{i\pi/4} {{z-i\lambda}\over {z-i\lambda^{-1}}} du  
\end{eqnarray}
 on the rhombic torus $$(2\cos \rho) w^{-2} = - (z - z^{-1} - 2i \sin \rho),$$ 
punctured at the two points where $z= i \lambda^{-1}$.

\begin{theorem}
For every $\rho \in (-\pi/2, \pi/2)$ and every 
$\lambda \in (0,1)\cup (1, \infty) $, the
Weierstrass data, (1.7)---defined on the rhombic torus
determined by $\rho$, with punctures where $z = i\lambda^{-1}$ on
the h-diagonal---produces in (1.2) a multivalued, regular,
complete minimal immersion of the punctured torus.

\begin{enumerate}
\item[(i)] Reflection in the v-diagonal (resp. h-diagonal) of
the torus induces an isometry of the minimal surface corresponding to
$180^\circ$ rotation about a vertical line (resp. horizontal line) in
${\Bbb R}^3$.


\item[(ii)] The immersion (1.2) is singly periodic if and only if 
$\lambda~\in~(0,1)$, $\rho~\in~(0,\pi/2)$ and

\end{enumerate}


\begin{eqnarray}
\Re \int_I dh & = & 0 \\
\Re \int_{II} i (g^{-1} + g) dh  & =  & 0,
\end{eqnarray}
where the paths $I$ and $II$ are as indicated in Figure~7. 
The period
\footnote{Since the
surface can be scaled by multiplying $dh$ by a nonzero real constant,
the value of $T$ has no geometric significance. For the representation we use
here, $T$ is given explicitly in (1.16) as a function of $\rho$ and $\lambda$.}
 is a vertical vector of the form $(0,0,T)$, 
where $T$ is given in (1.16) and
\begin{equation}
 T =\int_{``vert. diag.''} dh = \pm 2\pi i
(Residue_{z=i\lambda^{-1}}dh).
 \end{equation}

Conversely, any regular, complete minimal surface containing a vertical line, 
whose quotient by vertical translations has genus one, contains two horizontal
lines and has two helicoidal ends, is representable
with Weierstass data of the form
 (1.7) with $\lambda \in (0,1)$ and $\rho \in (0,\pi/2)$. 
\end{theorem}

\noindent {\bf Proof.} 
The "converse" part of the Theorem is
clear--- from the development in previous sections,
any such surface must be represented by the data
as given in (1.7) with punctures at the two points
 where $z = i\lambda^{-1}$. The statements that
 $\lambda \in (0,1)$ and 
that $\rho \in (0,\pi/2)$ will be proved below.

We begin by considering the one-form $dh$ in (1.7), which was
constructed to have its zeros precisely at the two points where $g = 0$ or
$g = \infty$ on the punctured torus. Hence the metric $$ ds = (\vert g
\vert + \vert g \vert^{-1}) \vert dh\vert$$ is regular. It is complete
because at the punctures, where $g = 0$ or $g = \infty$, $dh$ has a
simple pole.

Since the h- (resp. v-) diagonal makes an angle of $45^\circ$
(resp. -$45^\circ$) with the real axis, $\mu^*du = i\overline{du}$ (resp.
$\mu_{vert}^*du = -i\overline{du}$). From Section 1.2, we know that
$z\circ \mu_{vert}= z\circ \mu = -\overline{z}$. It follows 
from the definition of $dh$ in (1.7) that 
\begin{equation}
\mu^*dh = -\overline{dh} \  \hbox {\em  and } \  \mu_{vert}^*dh =
\overline{dh}.
\end{equation}

From Section 1.2 we have $w\circ \mu = -i \overline{w}$ (resp. 
 $w\circ \mu_{vert} = i \overline{w}$).  A computation using the
the definition of $g$ in (1.7) gives 
\begin{equation}
g\circ \mu = \overline{g} \  \hbox {  \em  and  } \  g\circ \mu_{vert} = 1/\overline{g}
.
\end{equation}

The Weierstrass representation (1.2) is $X(p)
  = X (p_0) + \Re \int_{p_0}^p \Phi $, 
$$
 \Phi :=\left (\phi_1, \phi_2, dh \right) =
 \left ( {1\over 2} (g^{-1} -g), ~{i\over 2} ( g^{-1} +g),
1\right)\,dh.
$$ 
From (1.11), (1.12) and the relation $ r\!_{_P} = \mu\circ \mu_{vert}$ 
we have: 
\begin{eqnarray}
\mu^*\Phi &=& \overline{\left(-\phi_1,\phi_2, -dh  \right)}\nonumber\\
\mu_{vert}^*\Phi &=& \overline{\left( -\phi_1,- \phi_2, dh \right)} \\
{r\!_{_P}}^*\Phi &=& \left(\phi_1,-\phi_2,-\phi_3 \right).\nonumber
\end{eqnarray}
Statement (i) of the theorem follows immediately from the Weierstrass
representation and these relations: writing $X =: \left(x_1,x_2,x_3 \right)$, 
and assuming, without loss of generality, that $X (p_0) = (0,0,0) $, then
\begin{eqnarray}
X\circ\mu(p) &= &\left(-x_1,x_2,-x_3 \right)(p)\nonumber\\
X\circ\mu_{vert}(p) &= & \left(-x_1,-x_2,x_3 \right)(p)\\
X\circ {r\!_{_P}}(p) &= & \left(x_1,-x_2,-x_3 \right)(p)\nonumber
\end{eqnarray}

We now address the period problem and begin by computing
the period at a puncture corresponding to an end.
  Because $g^{-1}dh$ and $ gdh$ both have 
their only pole at one of the punctures---and therefore
 no residue---$\phi_1$ and  $ \phi_2$ have no residues there.
Let $\alpha$ be a simple closed curve in the homotopy class of a puncture. (See
Figure~7.) Then 
$\int_\alpha \Phi = 2\pi i\left(0,0,Residue_{z=i\lambda^{-1}}dh \right).$
From (1.7) we have
\begin{eqnarray*}
dh= 
 e^{i\pi/4}{{z-i\lambda}\over{z-i\lambda^{-1}}}du  &=& 
 e^{i\pi/4}{{z-i\lambda}\over{z-i\lambda^{-1}}}({{z'}\over{z}})^{-1}
{{dz}\over{z}}\\&=&
e^{i\pi/4}{{z-i\lambda}\over{z}}
\left({{-2}\over{\cos\rho}}(z-z^{-1}-2i\sin\rho) \right)^{-1/2}
{{dz}\over{z-i\lambda^{-1}}}.\nonumber
\end{eqnarray*}
Hence 
\begin{eqnarray*}
Residue_{z=i\lambda^{-1}}dh &=& e^{i\pi/4} 
{{\lambda^{-1}-\lambda}\over{\lambda^{-1}}}
\left(
{{-4i}\over{\cos\rho}}({{\lambda +\lambda^{-1}}\over {2}}-\sin\rho) 
\right)^{-1/2}
\\ &=& 
i{-1}{{\sqrt{\cos\rho}}\over{2}}
 (1-\lambda^2)  
\left({{\lambda +\lambda^{-1}}\over{2}} -\sin\rho \right )^{-1/2}
\end{eqnarray*}
and 
\begin{equation}
Period_\alpha X :=\Re \int_\alpha  \Phi = (0, 0, \pm T),
\end{equation}
where
\begin{equation}
T:= \pi \sqrt{\cos\rho}(1-\lambda^2)  
\left({{\lambda +\lambda^{-1}}\over{2}} -\sin\rho \right )^{-1/2} \neq 0.
\end{equation}

Let $\alpha_1$and $\beta$ be as in
 Figure~7. Note that we may write
$\alpha = \alpha_1 -\mu\alpha_1$. We note for use below that because (from (1.13))  $\mu^* dh= -\overline{dh}$ and 
$\mu^*\phi_1= -\overline{\phi_1}$, 
$$ \Re \int_{\alpha_1} \phi_1 =\Re\int_{-\mu\alpha_1} \phi_1 \,\,\, \hbox{and}
\,\,\,   \Re\int_{\alpha_1} dh =\Re\int_{-\mu\alpha_1} dh.
$$
It follows from this and (1.15) that
\begin{equation}
\Re\int_{\alpha_1} \phi_1 =\Re\int_{-\mu\alpha_1} \phi_1 = 0 \,\,\,
\hbox{and}\,\,\,
\Re\int_{\alpha_1} dh =\Re\int_{-\mu\alpha_1} dh =\pm T/2.
\end{equation}

In Section~1.1 we chose, without loss of generality, to represent the surface 
so that the image of the center point of the torus was at the same height as
the two off-axis fixed points of the normal symmetry. Thus we must require
\begin{equation}
\Re \int_{\beta}dh = 0, 
\end{equation}
where $\beta$ is the curve in Figure~7.
By (1.14), it follows that $\Re \int_{\mu\beta}dh = 0$, so both off-axis
fixed points are at the same level. 

Any branch of $X$ maps the v-diagonal to a vertical line segment and,
since $g(O) = -1 $, the h-diagonal is mapped to a line parallel
to the $x_2-$axis. Hence the image under any branch of $X$ of the off-axis 
fixed points must lie on a line through $X(O)$, which is parallel to the
$x_1-$axis. That is
\begin{equation}
\Re \int_{\beta} \phi_2 = \Re\int_{\beta} ~{i\over 2} ( g^{-1} +g)dh = 0.
\end{equation}

Hence
$\Re \int_{\beta} \Phi =\left( \Re \int_{\beta} \phi_1,\, 0,\, 0 \right) $

For $\hat{\beta}:= {r\!_{_P}}\beta$, we have from (1.13) that
$\Re \int_{\hat {\beta}} \Phi =\left( \Re \int_{ \beta}\phi_1,\, 0,\, 0 \right)$.
 This means that on the closed cycle 
$B:=  +\beta -\hat \beta$,
\begin{equation}
Period_B X = \Re \int_{B} \Phi = (0,0,0).
\end{equation}

The four closed curves $\{B, \mu B, \alpha, \mu_{vert}\alpha\}$ form a
basis for the homology of the punctured torus. By (1.13) it follows
that $Period_B X =  (0,0,0)$ if and only if $Period_{\mu B} X =
(0,0,0)$.  This means that the {\em necessary} conditions (1.18) and
(1.19) are  {\em sufficient} to insure that X has no real period on
either $B$ or
 $\mu B$. By (1.15) and (1.16), X {\em always} has a nonzero vertical
period on $\alpha$ and---by symmetry, using (1.13)---on
$\mu_{vert}\alpha$. Therefore,  (1.18)
 and (1.19) are sufficient to guarantee that $X$ is singly periodic
with vertical period $(0,0,T)$.

Before translating (1.18) and (1.19) into (1.8) and (1.9), we will show
that (1.10) is valid. Observe that $B$ is homotopic to $D +
\gamma_1 + \alpha_1 +\gamma_2$. (See Figure~7.) We know
that $D$ is mapped into a vertical line line segment so
$\Re\int_{\gamma_1 + \alpha_1 +\gamma_2}\Phi$ is also a vertical vector.
The $\gamma_i$ lie on the horizontal  diagonal and are mapped  to horizontal
line segments.
 it follows from this and (1.17) that $$\Re\int_{D}\Phi = 
 - \left(0,0,\Re\int_{\alpha_1}dh  \right) = (0,0,-T/2)$$ This is
equivalent to (1.10) in the statement of Theorem~1.

Equations
(1.18) and (1.19) are the period conditions: necessary and sufficient conditions
for the surface in question to be singly periodic. We wish to express them
in a form more convenient for  subsequent calculation, i.e. as equations
(1.8) and (1.9). When $\lambda < 1 $, $\beta$ is homotopic to  
$F + I$ and to  $E + II $.
See Figure~7. Along $F$, we know
that $x_3$ is constant. Thus (1.18) is equivalent to 
$$
\Re \int_I dh = 0, 
$$
 which is (1.8).
Similarly, $x_2$ is constant along $E$, so (1.19) is
equivalent to
$$
Re \int_{II} {{i}\over{2}}(g^{-1} +g)dh = 0,
$$
which is (1.9). 

We now show that the period conditions cannot be satisfied when
 $\lambda~>~1.$
On $E$, $z= -it$, $0\leq t\leq 1$. We may use $z$
 to  parametrize $E$ by $t$:\, $ z\circ E(t) = -it$; $E(t) = f(t) e^{-i\pi/4}$, 
where $f(t)$ is some positive increasing function. Then 
\begin{eqnarray*}
dh (\dot E) &=& e^{i\pi/4} {{-it-i\lambda}\over{-it-i\lambda^{-1}}}e^{-i\pi/4}
\dot f(t) \nonumber \\
&=& {{t+\lambda}\over{t+\lambda^{-1}}}\dot f(t).
\end{eqnarray*}
Hence
$$
a := x_3(E(1)) = \Re \int_E dh > 0. 
$$

Along $II$, which begins at $b:=E(1)$, $z= -ie^{-it}$,  for $0\leq t\leq 
\pi/2 + \rho$. We may use $z$ to parametrize $II$ by $t$:
 $z\circ II(t) = -ie^{-it}$;  $II(t)= b - f(t) e^{+i\pi /4},$ where
$f(t)$ is, again, some positive increasing function. Then
\begin{eqnarray}
dh (\dot {II}) &=& e^{i\pi/4} {{-ie^{-it}-i\lambda}\over{-ie^{-it}-i\lambda^{-1}}}
e^{i\pi/4}(-\dot f(t))\nonumber\\
&=& -i\dot f(t){{e^{-it}+\lambda}\over{e^{-it}+\lambda^{-1}}}\\
&=& \left( (\lambda -\lambda^{-1} )
\sin t
+ i(\dots) \right){{\dot f(t)}\over{|e^{-it} + \lambda^{-1}|^2}}\nonumber
\end{eqnarray}
Thus $\Re dh(\dot{II})$ is strictly positive when $\lambda >1.$
This means that $x_3(t)$ is strictly increasing.
But since  $a = x_3(b)  > 0 $, this would imply that 
$\Re\int _{\beta}dh = \Re\int _{E +II}dh > 0$, violating (1.18).

%
%

This completes the proof of Theorem~1, with the exception of showing that
the period conditions cannot be solved for $\rho \leq 0$. This 
will be proved in Lemma~2 below.

\begin{remark} We note for use in Section 3 that (1.21) implies that when
$\lambda~<1$,  
$x_3\circ II(t)$ is  strictly decreasing.
\end{remark}

\subsection{The period conditions in terms of definite integrals}
In this section we  express the integrals in (1.8) and (1.9) 
in terms of explicit definite integrals. We will need this for the
existence proof in Section~2 as well as for the proof of 

\begin{lemma} The period condition  (1.9) is not solvable when
 $\rho \leq 0$.
\end{lemma}

Along paths I and II, where $z=e^{i\phi}$, we have
\begin{eqnarray*}
e^{i\pi/4}\,du&=&
e^{i\pi/4}\frac z{z'}\,\frac{dz}{z}=
e^{i\pi/4}\Bigl(-\frac{2}{
\cos\rho}(e^{i\phi}-e^{-i\phi}-2i\sin\rho)\Bigr)^{-1/2} i\,d\phi\cr
\noalign{\vskip3pt}&=&\cases{\displaystyle
\phantom{-i}\frac{\sqrt{\cos\rho}}{2}\frac{d\phi}{\sqrt{\sin\phi-\sin\rho}}
&for $\displaystyle\rho<\phi<{\pi/2}$,\cr\noalign{\vskip3pt}\displaystyle
{-i}\frac{\sqrt{\cos\rho}}{2}\frac{d\phi}{\sqrt{\sin\rho-\sin\phi}}
&for $\displaystyle-{\pi/2}<\phi<\rho$}
\end{eqnarray*}
and
\begin{eqnarray*}
\frac{z-i \lambda}{ z-i \lambda^{-1}}&=& \frac{e^{i\phi}-i \lambda}{
e^{i\phi}-i \lambda^{-1}}\, \frac{\lambda e^{-i\phi}+i}{\lambda
e^{-i\phi}+i}\cr\noalign{\vskip3pt} &=&
\lambda\,\frac{2-(\lambda+\lambda^{-1})
\sin\phi+i(\lambda^{-1}-\lambda)\cos\phi}{\lambda+\lambda^{-1}-2\sin\phi}.
\end{eqnarray*}
Let $\Lambda=\lambda+\lambda^{-1}$.  Using (1.7) and the expressions
above, we see that the first period condition (1.8) is equivalent to
\begin{eqnarray}
\int_\rho^{\frac\pi2}\frac{2- \Lambda\sin \phi}{ \Lambda-2\sin \phi}\,
\frac{d\phi}{\sqrt{\sin\phi-\sin\rho}}=0.
\end{eqnarray}

To write the second period condition (1.9), we first compute
\begin{eqnarray*}
g+\frac 1g &=& \frac{w-re^{i\pi/4}}{w+re^{i\pi/4}}+
\frac{w+re^{i\pi/4}}{w-re^{i\pi/4}}= 2\,\frac{w^2+ir^2}{w^2-ir^2}\\[4pt]&=&
2\,\frac{\Lambda-4\sin\rho+2\sin\phi}{\Lambda-2\sin\phi}.
\end{eqnarray*}
The last equality follows from (1.6) and the
definition of $\Lambda=\lambda+\lambda^{-1}$.
  Together with the expression for $dh$ computed 
above, this shows that condition (1.9) is equivalent to
\begin{eqnarray}
\int_{-\frac\pi2}^{\rho}
\frac{\Lambda-4\sin\rho+2\sin\phi}
{\Lambda-2\sin\phi}\, \frac{2- \Lambda\sin \phi}{ \Lambda-2\sin \phi}\,
\frac{d\phi}{\sqrt{\sin\rho-\sin\phi}}=0.
\end{eqnarray}

For $\rho\leq0$, no factor of the integrand of (1.23) changes
sign, so the second period condition is not solvable for
$-\pi/2<\rho\le0$.


\def\tfrac#1#2{{\textstyle\frac{#1}{#2}}}

\section{Existence}
\renewcommand{\theequation}{2.\arabic{equation}}
\setcounter{equation}{0}
\begin{theorem}
There exists $(\rho_0,\lambda_0)\in(0,\,\pi/2)\times(0,1)$ satisfying
the period conditions (1.8) and (1.9). 
\end{theorem}

{\bf Proof.} In Section~1.6, we wrote  the period conditions (1.8) and
(1.9) as the definite integrals (1.22) and (1.23). We will work
with these relations without reference to the geometry of their
derivation:

\begin{eqnarray}
\int_\rho^{\tfrac\pi2}\frac{2- \Lambda\sin \phi}{ \Lambda-2\sin \phi}\,
\frac{d\phi}{\sqrt{\sin\phi-\sin\rho}}&=&0;
\end{eqnarray}

\begin{eqnarray}
\int_{-\tfrac\pi2}^{\rho}
\frac{\Lambda-4\sin\rho+2\sin\phi}
{\Lambda-2\sin\phi}\, \frac{2- \Lambda\sin \phi}{ \Lambda-2\sin \phi}\,
\frac{d\phi}{\sqrt{\sin\rho-\sin\phi}}&=&0,
\end{eqnarray} 
where $\Lambda =\lambda+\lambda^{-1}$. We need to find a pair $(\rho,~
\lambda)$ solving both of these equations simultaneously. 

The  proof goes as follows.  
Using the intermediate value theorem, we  show in Section~2.1 
that condition (2.1) can be
solved for $(\rho,\lambda(\rho))$, where $\lambda(\rho)$ is a
differentiable function of $\rho$ and $0<\rho<\pi/2$.  Then we show
in Section~2.2
that the period integral in (2.2) changes sign along the graph of
$\lambda=\lambda(\rho)$. By the intermediate value theorem again, there
exists $(\rho_0,\lambda(\rho_0))$ at which both (2.1) and (2.2) hold.
\qed

\subsection{Solution of the  first period integral as a function of $\rho$}

Consider, for $0<\rho<\pi/2$, the first period integral in (2.1).
Observe  that the differentiable function
$$
F(\rho, \Lambda)=\int_\rho^{\tfrac\pi2}
\frac{2-\Lambda \sin\phi}{\Lambda-2 \sin\phi}\,
\frac{d\phi}{\sqrt{\sin\phi-\sin\rho}},
$$
satisfies
\begin{eqnarray}
F(\rho,2)>0\quad\hbox{and}\quad
F\Bigl(\rho,\frac2{\sin\rho}\Bigr)<0,
\end{eqnarray}
and that $F(\rho,\Lambda)$ is strictly decreasing in $\Lambda$.  Therefore
we may define a differentiable function $\Lambda(\rho)$, for
$0<\rho<\pi/2$, by the condition
$$
F(\rho,\Lambda(\rho))=0.
$$
By (2.3) we have
\begin{eqnarray}
2<\Lambda(\rho)<\frac{2}{\sin\rho}.
\end{eqnarray}

In order to estimate the second period integral along the graph of
$\Lambda=\Lambda(\rho)$, we need to control $\Lambda(\rho)$.

\begin{lemma}
For $\rho\in(0,\pi/2)$ we have
$$
2<\Lambda(\rho)<\min\Bigl(\frac{2}{\sin\rho},\,8\Bigr).
$$
For $\rho\in(\pi/2-\epsilon,\,\pi/2)$,
where $\epsilon$ is sufficiently small, we have
$$
2<\Lambda(\rho)<2+(1-\sin\rho).
$$
\end{lemma}

\noindent {\bf Proof.} 
Since we already have (2.3) and (2.4), it is enough to show that
$F(\rho,8)<0$ for $\rho\in(0,\pi/2)$, and
$F(\rho,\,2+(1-\sin\rho))<0$ for $\rho$ near $\pi/2$.
Define $\phi_\Lambda $ to be the zero of the integrand of $F$, so that
$$
\sin \phi_\Lambda=\frac2\Lambda.
$$
Then 

\begin{eqnarray}
F(\rho,\Lambda) &=& \int_\rho^{\phi_\Lambda} \frac{2-\Lambda
\sin\phi}{\Lambda-2 \sin\phi}\,
\frac{d\phi}{\sqrt{\sin\phi-\sin\rho}}\nonumber\\ &-&
\int_{\phi_\Lambda}^{\tfrac\pi2} \frac{\Lambda \sin\phi-2}{\Lambda-2\sin\phi}\,
 \frac{d\phi}{\sqrt{\sin\phi-\sin\rho}}.
\end{eqnarray}

Both integrands are positive. Since $(2-\Lambda \sin\phi)/(\Lambda-2
\sin\phi)$ is decreasing in the interval $(\rho,\phi_\Lambda)$ and
$(\cos\phi)/(\cos\phi_\Lambda)\ge1$, we have
\begin{eqnarray*}
\int_\rho^{\phi_\Lambda} \frac{2-\Lambda \sin\phi}{\Lambda-2
\sin\phi}\, \frac{d\phi}{\sqrt{\sin\phi-\sin\rho}}&\le&
\frac{2-\Lambda \sin\rho}{\Lambda-2 \sin\rho}
\int_\rho^{\phi_\Lambda}\frac {\cos\phi\,/\cos\phi_\Lambda}
{\sqrt{\sin\phi-\sin\rho}}\,d\phi\cr\noalign{\vskip5pt} &=&
\frac{2-\Lambda \sin\rho}{\Lambda-2 \sin\rho}\,
\frac{2\sqrt{\sin\phi_\Lambda-\sin\rho}}
{\cos\phi_\Lambda}\cr\noalign{\vskip5pt}
&=&\frac{2-\Lambda \sin\rho}{\Lambda-2 \sin\rho}\cdot2
\sqrt{\frac{2/\Lambda-\sin\rho}{1-4/\Lambda^2}}.
\end{eqnarray*}

Turning our attention to the second integral in (2.5), we first estimate
a term in the integrand. Let
 $$
f(\phi):=\frac{\Lambda \sin \phi-2}{\Lambda-2\sin\phi}. 
 $$ 
By the
definition of $\phi_\Lambda$, $f(\phi_\Lambda)=0$ and clearly
$f(\pi/2)=1$.
 It is easy to check that $f'(\phi)>0$  on $(\phi_\Lambda,\,\pi/2)$, so
$f(\phi)$ is increasing on this interval. Also one can compute
that $f'(\phi_\Lambda)= 1/\cos(\phi_\Lambda) >1/(\pi/2 -\phi_\Lambda)$.
Furthermore,
  $f''(\phi_\Lambda)>0$ and $f''(\phi)$ has only one zero on
$(\phi_\Lambda,\,\pi/2)$. The linear function $$
l(\phi):=\frac{\phi-\phi_\Lambda}{\pi/2-\phi_\Lambda} $$ satisfies
$l(\phi_\Lambda)=0=f(\phi_\Lambda)$ and  $l(\pi/2)=1=f(\pi/2)$.  Since
$l'(\phi) = 1/(\pi/2-\phi_\Lambda)$, we have
$f'(\phi_\Lambda)>l'(\phi_\Lambda)$. Because $l''(\phi) = 0$, 
$f(\phi)>l(\phi)$ at least up to the first zero of $f''$. On the
remaining subinterval, $f$ is concave, so it is above its secant, which
is in turn above $l$. It follows that 
$$
f(\phi)=\frac{\Lambda \sin \phi-2}{\Lambda-2\sin\phi}\ge\frac{
\phi-\phi_\Lambda}{\pi/2-\phi_\Lambda}
$$
on $(\phi_\Lambda,\pi/2)$.

We now can now estimate the integral in (2.5).
\begin{eqnarray*}
\int_{\phi_\Lambda}^{\tfrac{\pi}{2}}
\frac{\Lambda \sin \phi-2}{\Lambda-2\sin\phi}
\,\frac{d\phi}{\sqrt{\sin\phi-\sin\rho}}&\ge&
\int_{\phi_\Lambda}^{\tfrac{\pi}{2}}
\frac{\phi-\phi_\Lambda}{\pi/2-\phi_\Lambda}\,\frac{d\phi}{\sqrt{1-\sin\rho}}
\\[4pt]&=&\frac12\,\frac{\pi/2-\phi_\Lambda}{\sqrt{1-\sin\rho}}\ge
\frac12\,\frac{\sin(\pi/2-\phi_\Lambda)}{\sqrt{1-\sin\rho}}\\[4pt]
&=&\frac12\,\frac{\cos\phi_\Lambda}{\sqrt{1-\sin\rho}}=
\frac12\,\frac{\sqrt{1-4/\Lambda^2}}{\sqrt{1-\sin\rho}}.
\end{eqnarray*}
These two integral estimates and (2.5) give
$$
F(\rho,\Lambda)\le2\,\frac
{2-\Lambda \sin \rho}{\Lambda-2\sin\rho}\sqrt{\frac
{2/\Lambda-\sin\rho}{1-4/\Lambda^2}}-\frac12\sqrt{\frac{1-4/\Lambda^2}
{1-\sin\rho}}.
$$
Notice that the right-hand side is a decreasing function of $\rho$ for
$\rho\in (0,\,\pi/2)$.  At $\Lambda=8$ we have
$$
F(\rho,8)\le F(0,8)\le 2\cdot\frac 28\cdot\sqrt{\frac{2/8}{1-4/8^2}}-
\frac12\sqrt{1-\frac{4}{8^2}}<0.
$$
This proves the first part of Lemma~3, that is,
$$
2 < \Lambda(\rho)<\min\Bigl(\frac2{\sin\rho} \, , \,8\Bigr)
$$
for $\rho\in(0,\,\pi/2)$.

\bigskip

We now wish to estimate $\Lambda(\rho)$ for $\rho$ near $\pi/2$.  Let
$\Lambda_k:=2+k(1-\sin\rho)$.  Then
$$
F(\rho,\Lambda_k)\le
2\,\frac{(2-k\sin\rho)^{3/2}}{2+k}\sqrt{\frac{\Lambda_k}{k(4+k(1-\sin\rho))}}
 -\frac{1}{2 \Lambda_k }\sqrt{k(4+k(1-\sin\rho))}.
$$
As $\rho\to\pi/2$, we have $\Lambda\to 2$ by (2.4), and
$$
\limsup_{\rho\to\pi/2} F(\rho,\Lambda_k)\le\frac{(2-k)^{3/2}}{2+k}
\sqrt{\frac2k}-\frac{\sqrt k}{2}.
$$
If $k=1$, then
$$
\limsup_{\rho\to\pi/2} F(\rho,\Lambda_k)\le\frac{\sqrt2}3-\frac{1}2<0.
$$
Together with the fact that $F(\rho,2)>0$, this shows that there
exists $\epsilon>0$ such that
$$
2<\Lambda(\rho)<2+(1-\sin\rho)
$$
for $\rho\in(\pi/2-\epsilon,\,\pi/2)$, which is the second estimate of
the lemma.\qed

\subsection{Solution of the second period integral along the graph
$\Lambda=\Lambda(\rho)$} 

We now show that the second period integral in (2.2),
$$
G(\rho,\Lambda)=\int_{-\tfrac\pi2}^{\rho} \frac{\Lambda-4\sin\rho+2\sin\phi}
{\Lambda-2\sin\phi}\, \frac{2- \Lambda\sin \phi}{ \Lambda-2\sin \phi}\,
\frac{d\phi}{\sqrt{\sin\rho-\sin\phi}},
$$
changes sign along $\Lambda=\Lambda(\rho)$.  This will guarantee that
there exists $(\rho_0,\,\Lambda(\rho_0))$ at which both the period
conditions (2.1) and (2.2), or, equivalently, (1.8) and (1.9), are
satisfied.

Since $2<\Lambda<8$ by Lemma~3, we see that the integrand of
$G(\rho,\Lambda)$ converges to a positive function as $\rho$ decreases to zero.
Therefore $G(\rho,\Lambda)>0$ for $\rho$ near 0. We will now show that
$G(\rho,\Lambda)<0$ for $\rho$ near $\pi/2$.

As $\rho\to\pi/2$, the function $G(\rho,\Lambda)$ may or may not have
a limit, depending on how fast $\Lambda(\rho)$ converges to 2.  So we
estimate $\limsup_{\rho\to\pi/2}G(\rho,\Lambda)$ instead, using the
second part of the estimates in Lemma~3.

For each $\rho$ smaller than and sufficiently close to $\pi/2$, define
$\phi_\rho\in(0,\rho)$ by the condition $$
\Lambda-4\sin\rho+2\sin\phi_\rho=0.  $$By Lemma~3, all the terms in the
integrand of  $G(\rho,\Lambda)$ are positive on $(-\pi/2, \rho)$, except
$\Lambda-4\sin\rho+ 2\sin\phi$, whose sign is the sign of the
integrand.  This allows us to estimate $G(\rho,\Lambda)$ by two integrals
we can control:  $$ G(\rho,\Lambda)=\int_{-\tfrac\pi2}^{\phi_\rho}+
\int_{\phi_\rho}^\rho\le \int_{-\tfrac\pi2}^0+ \int_{\phi_\rho}^\rho,
$$ since the integrand of $G(\rho,\Lambda)$ is negative on
$(-\pi/2,\phi_\rho) \supset (0,\phi_\rho)$.

The integral on $(-\pi/2,\,0)$ is easy to control.  It is continous at
$\rho=\pi/2$ and $\Lambda=2$, and we have
\begin{eqnarray}
&&\lim_{\rho\to\pi/2}
\int_{-\tfrac{\pi}{2}}^{0}\frac{\Lambda-4\sin\rho+2\sin\phi}{\Lambda-2\sin\phi}
\,\frac{2-\Lambda\sin\phi}{\Lambda-2\sin\phi}\,\frac{d\phi}{\sqrt{\sin\rho-\sin\phi}}
\cr
&&\qquad=
 -\int_{-\tfrac{\pi}{2}}^{0}\frac{d\phi}{\sqrt{1-\sin\phi}}\cr &&\qquad\le
 -\int_{-\tfrac{\pi}{2}}^{0}\frac{d\phi}{\sqrt{1-\phi}} =
2(1-\sqrt{1+\pi/2})\approx -1.2067.
\end{eqnarray}

On the interval $(\phi_\rho,\rho)$, the function $\sin\phi$ is
monotonic increasing.  It is easy to check that $$
g(\sin\phi):=\frac{\Lambda-4\sin\rho+2\sin\phi}{\Lambda-2\sin\phi} $$
is a convex increasing function of $\sin\phi$ with range $[0,1]$.
Therefore $$ g(\sin\phi)\le
g(\sin\phi_\rho)+\frac{g(\sin\rho)-g(\sin\phi_\rho)}{\sin\rho
 -\sin\phi_\rho}(\sin\phi-\sin\phi_\rho).  $$ That is, $$
\frac{\Lambda-4\sin\rho+2\sin\phi}{\Lambda-2\sin\phi}\le
\frac{\sin\phi-\sin\phi_\rho}{\sin\rho -\sin\phi_\rho}.  $$ Since
$\Lambda>2$ by (2.4), we have $$
\frac{2-\Lambda\sin\phi}{\Lambda-2\sin\phi}\le 1, $$ so we estimate the
integral on $(\phi_\rho,\rho)$ as follows:  \begin{eqnarray*}
&&\int_{\phi_\rho}^\rho
\frac{\Lambda-4\sin\rho+2\sin\phi}{\Lambda-2\sin\phi}
\,\frac{2-\Lambda\sin\phi}{\Lambda-2\sin\phi}\,\frac{d\phi}{\sqrt{\sin\rho-\sin\phi}}\\[4pt]
&&\qquad\le \frac1{\sin\rho-\sin\phi_\rho} \int_{\phi_\rho}^
\rho\frac{\sin\phi-\sin\phi_\rho}{\sqrt{\sin\rho-\sin\phi}}\,d\phi\\[4pt]
&&\qquad\le \frac1{\sin\rho-\sin\phi_\rho} \int_{\phi_\rho}^
\rho\frac{\sin\phi-\sin\phi_\rho}{\sqrt{\sin\rho-\sin\phi}}\,
\frac{\cos\phi}{\cos\rho}\,d\phi\\[4pt] &&\qquad=
\frac1{(\sin\rho-\sin\phi_\rho)\cos\rho}\cdot
\frac43\,(\sin\rho-\sin\phi_\rho)^{3/2}
=\frac43\frac{\sqrt{\sin\rho-\sin\phi_\rho}}{\sqrt{1-\sin^2\rho}}.\cr
\end{eqnarray*}  
Since $\Lambda+4\sin\rho+2\sin\phi_\rho = 0$, by the
definition of $\phi_\rho$, and
 $\Lambda\le 2+(1-\sin\rho)$ by Lemma~3, we know that 
$$
2(\sin\rho-\sin\phi_\rho)=\Lambda-2\sin\rho\le3(1-\sin\rho).$$
Therefore
$$
\frac43
\sqrt{\frac{\sin\rho-\sin\phi_\rho}{1-\sin^2\rho}}\le
\frac43
\sqrt{\frac{\frac32(1-\sin\rho)}{1-\sin^2\rho}}=\frac
{4\sqrt{\frac32}}{3\sqrt{1+\sin\rho}}.
$$
Thus we have
$$
\limsup_{\rho\to\pi/2}
\int_{\phi_\rho}^{\rho}\frac{\Lambda-4\sin\rho+2\sin\phi}{\Lambda-2\sin\phi}
\,\frac{2-\Lambda\sin\phi}{\Lambda-2\sin\phi}\,\frac{d\phi}{\sqrt{\sin\rho-\sin\phi}}
\le\frac{4\sqrt{\frac32}}{3\sqrt2}\approx 1.1547.
$$
 This estimate,
together with (2.6), shows that $$ \limsup_{\rho\to\pi/2}
G(\rho,\Lambda)<0.  $$ Therefore $G(\rho,\Lambda)<0$ for $\rho$ near
$\phi/2$, as desired.

\section{Embeddedness}
\renewcommand{\theequation}{3.\arabic{equation}}
\setcounter{equation}{0}
\begin{theorem}
Any minimal immersion with Weierstrass data (1.7),
which is singly periodic (i.e. satisfies (1.8) and (1.9)), must be an embedding.
\end{theorem}

 
We begin by giving a sketch of the proof, which we hope is easy to understand
in outline, with references to subsequent sections where the details
are carried out. The key idea is to cut a translational fundamental domain of the surface
into four congruent pieces and then show that each piece is a graph over
a domain in the $(x_1, x_2)$-plane. The graphs meet only along their 
embedded boundary arcs. 

To produce a translational fundamental domain, we cut the torus on
the h-diagonal from one end to the other, along the part of the
diagonal that contains the vertex (labelled $O'$ in Figure~8).  The
Weierstrass integral (1.2) is single valued on this slit torus and
produces as image a fundamental domain bounded above and below by
horizontal straight lines that cross the vertical axis of the surface.
These lines are the images of the slit. 

We now want to cut up the slit torus into regions whose images are
graphs over the $(x_1, x_2)$-plane. This means that the cuts must
include all points where the Gauss map is horizontal, i.e. all points
where $|g|=1$. This set consists of the v-diagonal, together
with the line (on the torus) parallel to the h-diagonal and
passing through the branch points of $g$ (labelled $b, \hat b$ in
Figure~8) on the v-diagonal. This decomposes the slit torus into
two rectangular pieces. We cut once more along the h-diagonal
from end to end along the segment through the center point (labelled
$O$ in Figure~8), producing four rectangles.

From Theorem~1i) in Section~1, we know that the images of the four
rectangles under the Weierstrass integral $X$ in (1.2) are congruent.
On the interior of each rectangle, $|g|\neq 1$, so the projection onto
the $(x_1, x_2)$-plane of its image under $X$ is an immersion. We will
prove in Proposition~1 below that, in fact, the projection is
one-to-one, i.e. the image of the interior of each rectangle is a
graph.  This requires us to control carefully the behavior of the image
of the boundaries of the rectangles. These boundaries consist of three different
parts (See Figure~9): 
 \begin {itemize} \item The image of half the
h-diagonal. This consists of two parallel horizontal
half-rays. One, the image of $H_1$, begins at $X(O)$, which we will
place at the origin in $R^3$.  The other, the image of $H_2$, begins at
$X(O')$ on the $x_3$-axis; \item The image of two segments of the
v-diagonal. (These are $E$ and $\hat{E}$, which add up to half
of the v-diagonal.  Each piece connects $X(b)$ or $X(\hat{b})$ to a point
of intersection of the v- and h-diagonals.) This
consists of two disjoint segments on the vertical axis, separated by
half a period; \item The image of the remaining arc, on which $|g| = 1
$ (labelled $C$ in Figure~9). This curve joins the branch points of $g$
on the vertical axis. The third component of the curve is a decreasing
function on $C$, which implies that the curve is embedded.  Since the
arc contains no branch points of $g$ in its interior, it projects to a
curve in the $(x_1, x_2)$-plane whose tangent vector rotates at a speed
that is never zero. This means that the projected loop, $c$,  is
everywhere locally convex. Since the total turning of its tangent
vector is less than $2\pi$, it is convex.  \end{itemize}

Putting together the above information, which will be established in
Proposition 1, with some more details and arguments, we will show that the
boundary of each rectangle is embedded by $X$ and that the interior of each
rectangle is embedded as a graph over a halfplane,  minus the interior of
the convex loop $c$.\footnote{Note that if the normal to a minimal surface
is horizontal  along a curve that projects to a strictly convex curve, 
$c$, the minimal surface near the curve projects  to the outside of
$c$.} The curve $c$ has one point in common with the boundary
of the halfplane.(See Figure~9, bottom left.)

The proof of Theorem 3 is completed in Section~3.3, where it is shown
that the four open graphs are disjoint and that the union of their boundary
curves is embedded.

\subsection{The statement of Proposition 1} As described in the
previous section, we cut the torus into four rectangles by removing the
points where $|g| = 1$,  as well as the h-diagonal, along
which $g$ is real.
 The function $$g = {{w-r e^{i\pi/4}}\over {w + re^{i\pi/4}}}$$ is
 unitary precisely on the set of points where $w$ takes values in
$e^{-i\pi/4} {\Bbb R}$, and we know from Section~1.3 that this set
consists of the v-diagonal, together with the line orthogonal to
this diagonal passing through the branch points of $g$.  (Note that
this curve passes through the off-axis fixed points of the normal
symmetry, which are also branch points of $z$.)   Reflection in the
diagonals and rotation about $O\in T^2$  induce Euclidean motions that
act transitively on the four images of these regions under $X$. (See
Figures~8~and~9).

\begfig
\hspace{0.5cm}
\epsfxsize=1.2in 
\epsffile{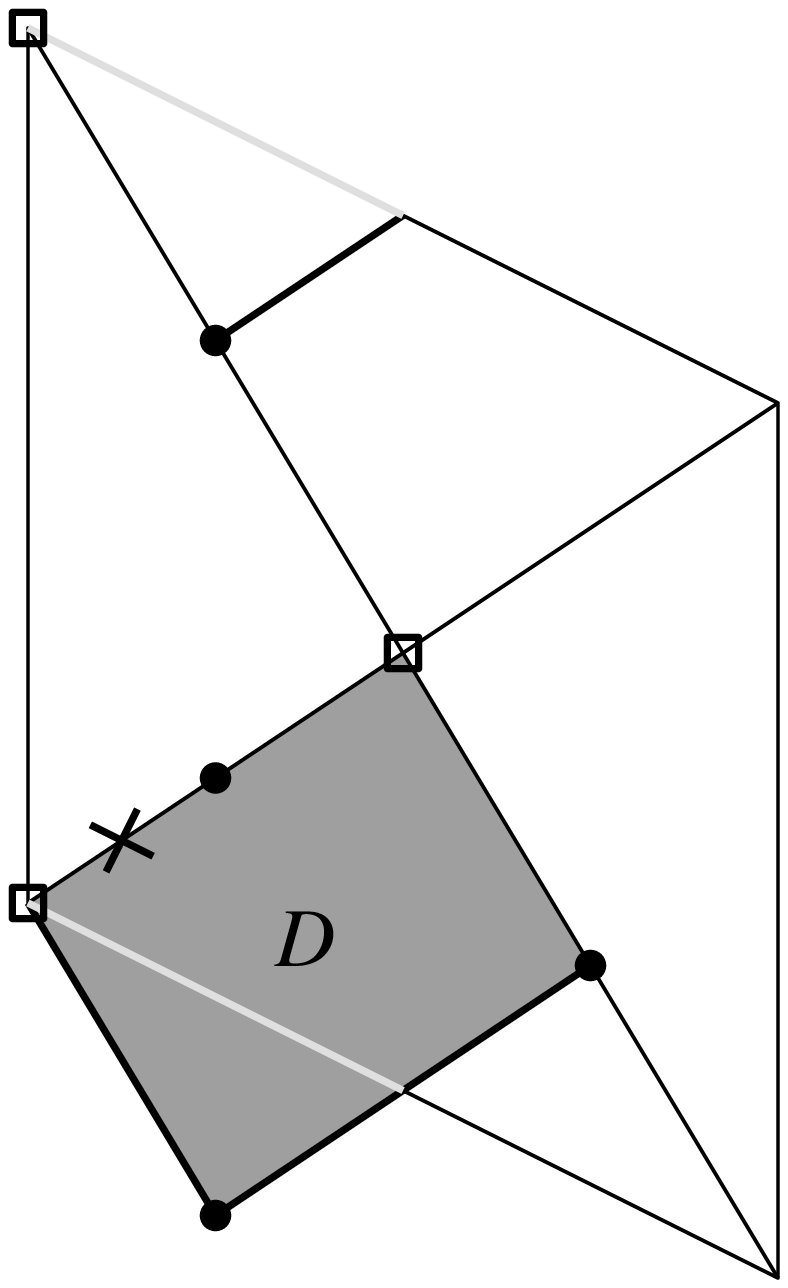}
\hspace{0.4cm}
\epsfxsize=1.5in 
\epsffile{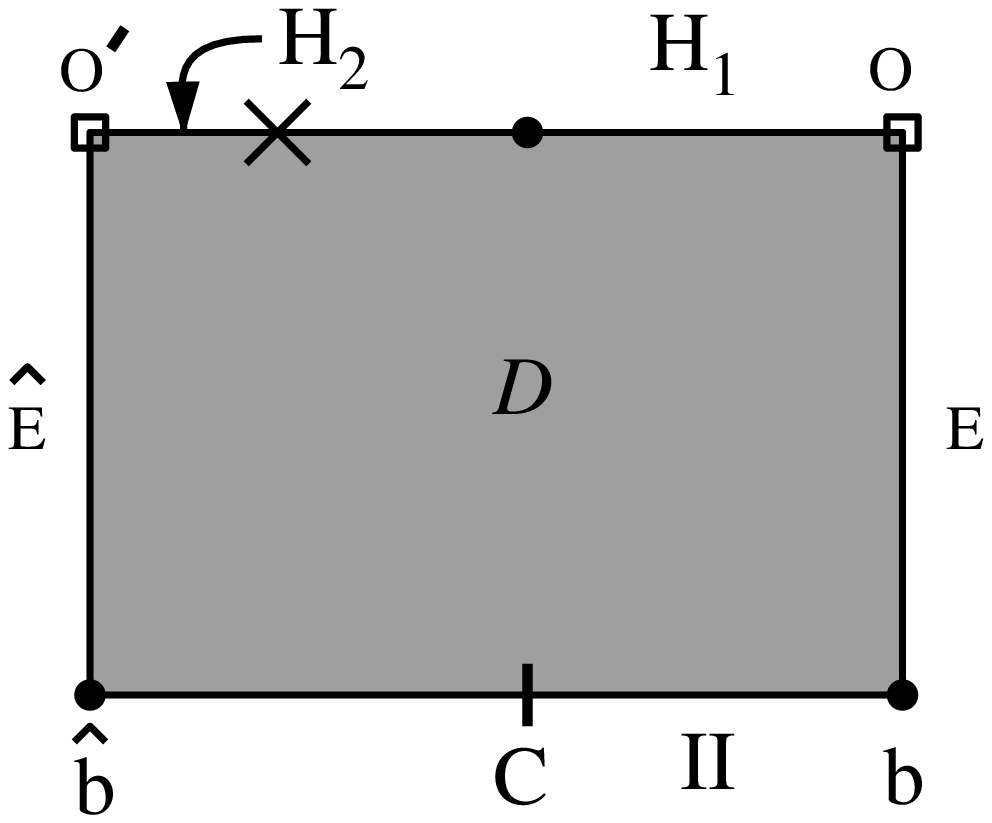}
\hspace{0.4cm}
\epsfxsize=1.2in 
\epsffile{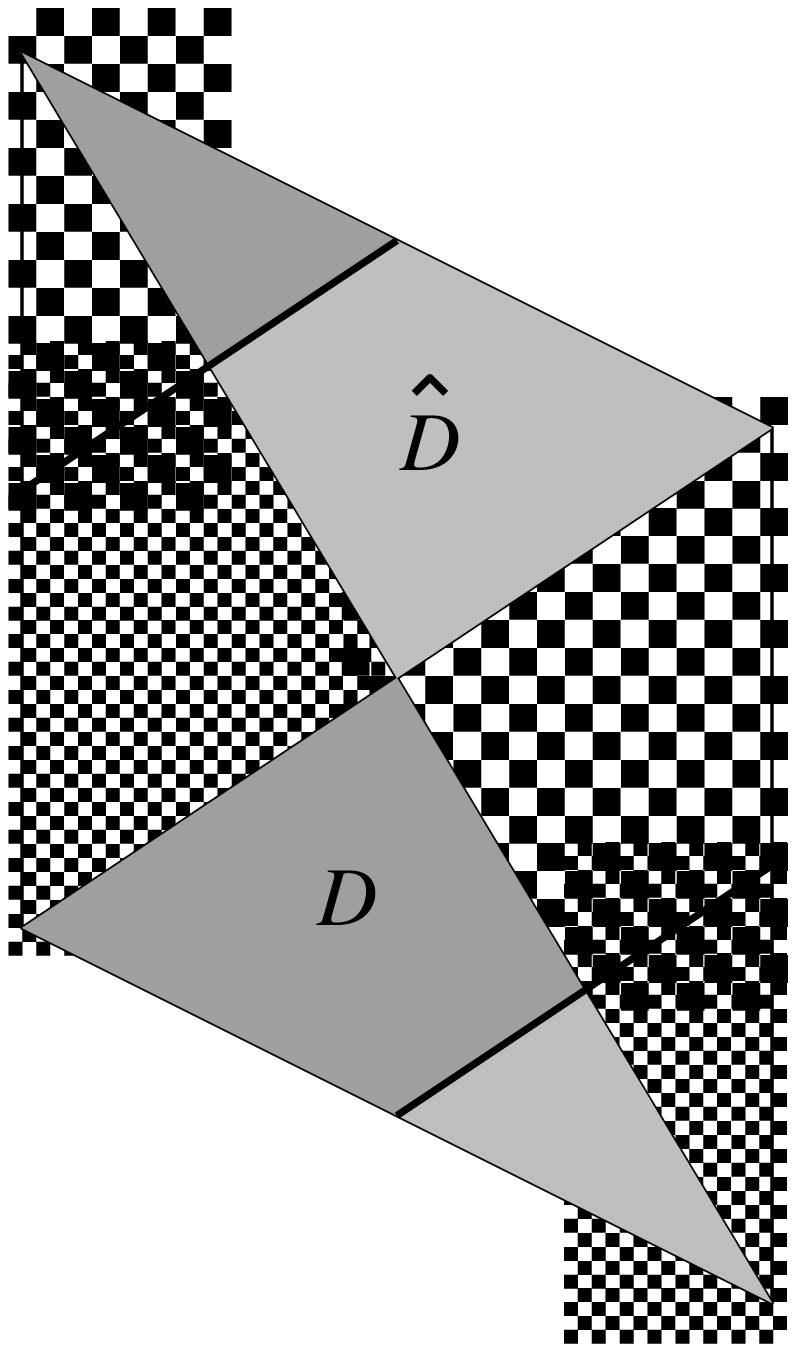}
\hfill
\vspace{0.2cm}

\def\thecaption{
{\em Left and Center}: The region $\mathbf{\cal D}$, whose image is a
graph over a halfplane minus a compact convex set.  {\em Right}: The
four regions, each congruent to  $\mathbf{\cal D}$ by automorphisms
that extend to Euclidean symmetries. The region ${\mathbf{\hat{\cal
D}}}= r_P \mathbf{\cal D}$ is referred to in the proof of Theorem~3. 
The puncture (x) divides the segment $OO'$ into the subsegments $H_1$
and $H_2$.}
\caption{\thecaption}
\endfig

Each of the four regions has the property that  $\vert g\vert \neq1$ 
on its interior.  We will work with the region labeled
$\mathbf{\cal D}$ in Figure~8.  Its Gaussian image satisfies $\vert g\vert >1$ on
the interior, with $\vert g\vert =1$ on $E  \cup \hat E \cup C$, and
 $g(H_1\cup H_2)\subset{\Bbb R} $.

\begfig
\vspace{-.5cm}
\hspace{-0.2cm}
\epsfxsize=2.4in 
\epsffile{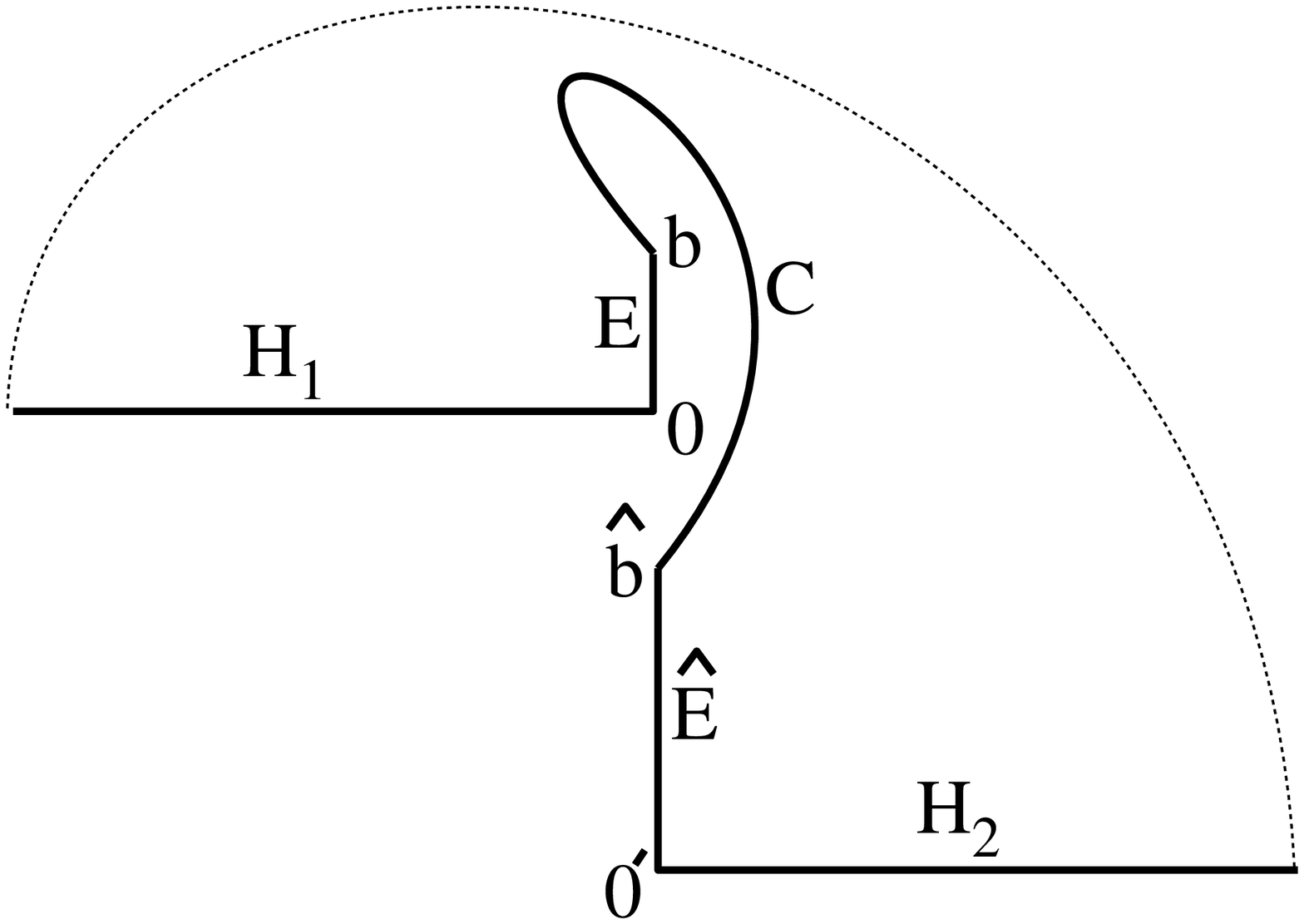}
\hspace{0.4cm}
\epsfxsize=2.4in 
\epsffile{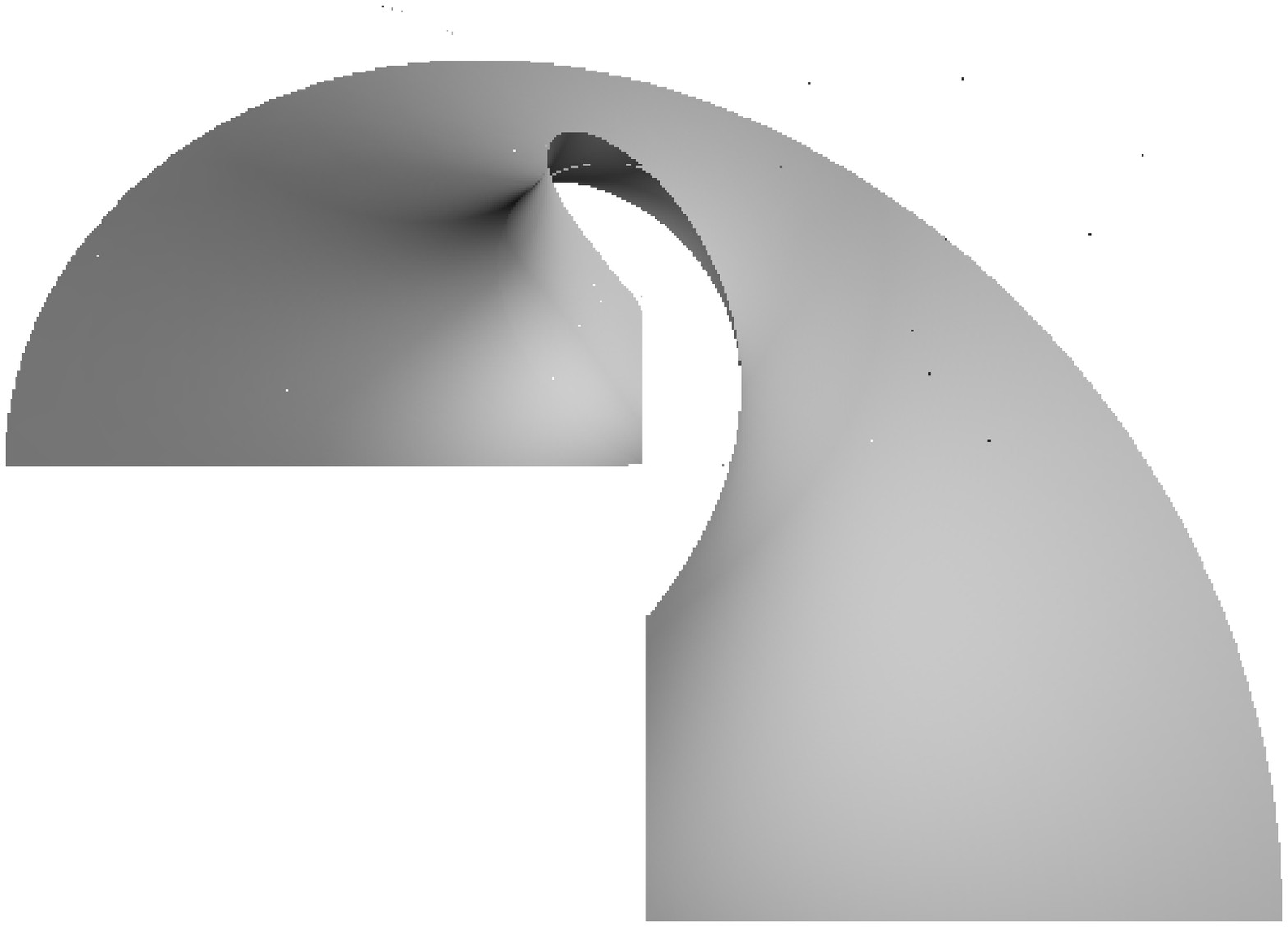}
\hfill
\vspace{-0.2cm}

\hspace{-0.2cm}
\epsfxsize=2.4in 
\epsffile{./figures/fig14.eps}
\hspace{0.4cm}
\epsfxsize=2.4in 
\epsffile{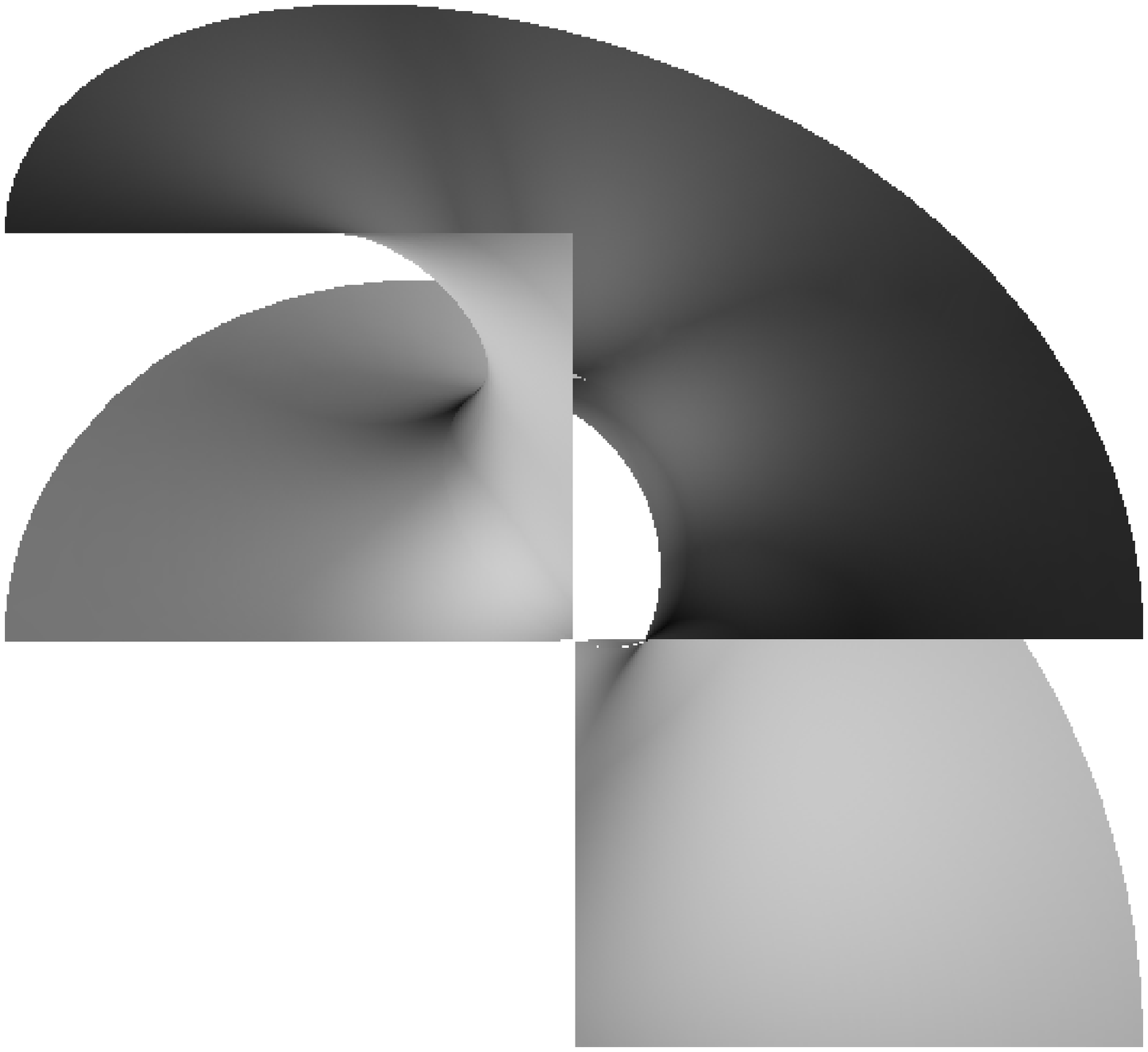}
\hfill
\vspace{-0.2cm}

\hspace{-0.2cm} \epsfxsize=2.4in \epsffile{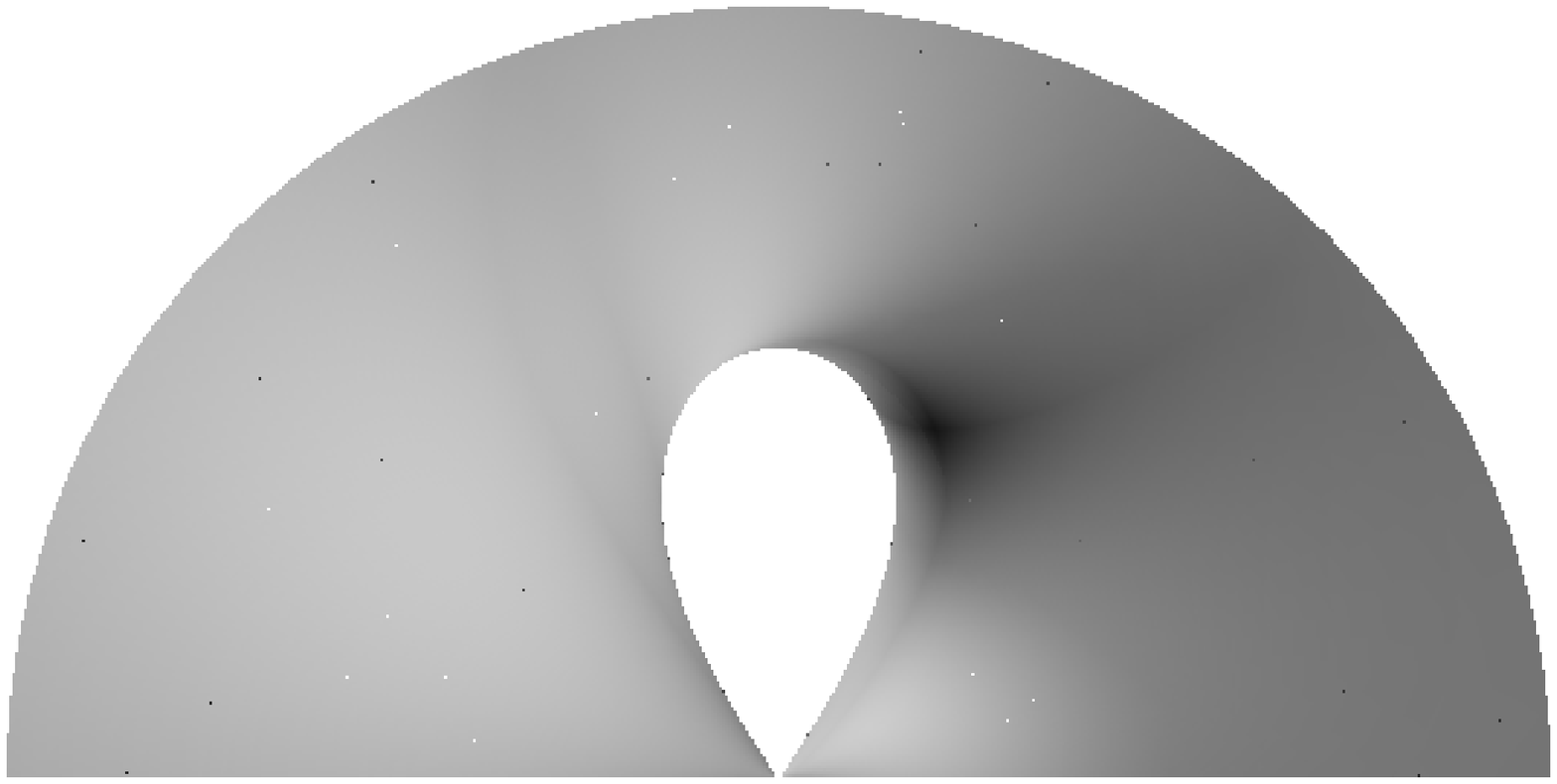}
\hspace{0.4cm} \epsfxsize=2.4in \epsffile{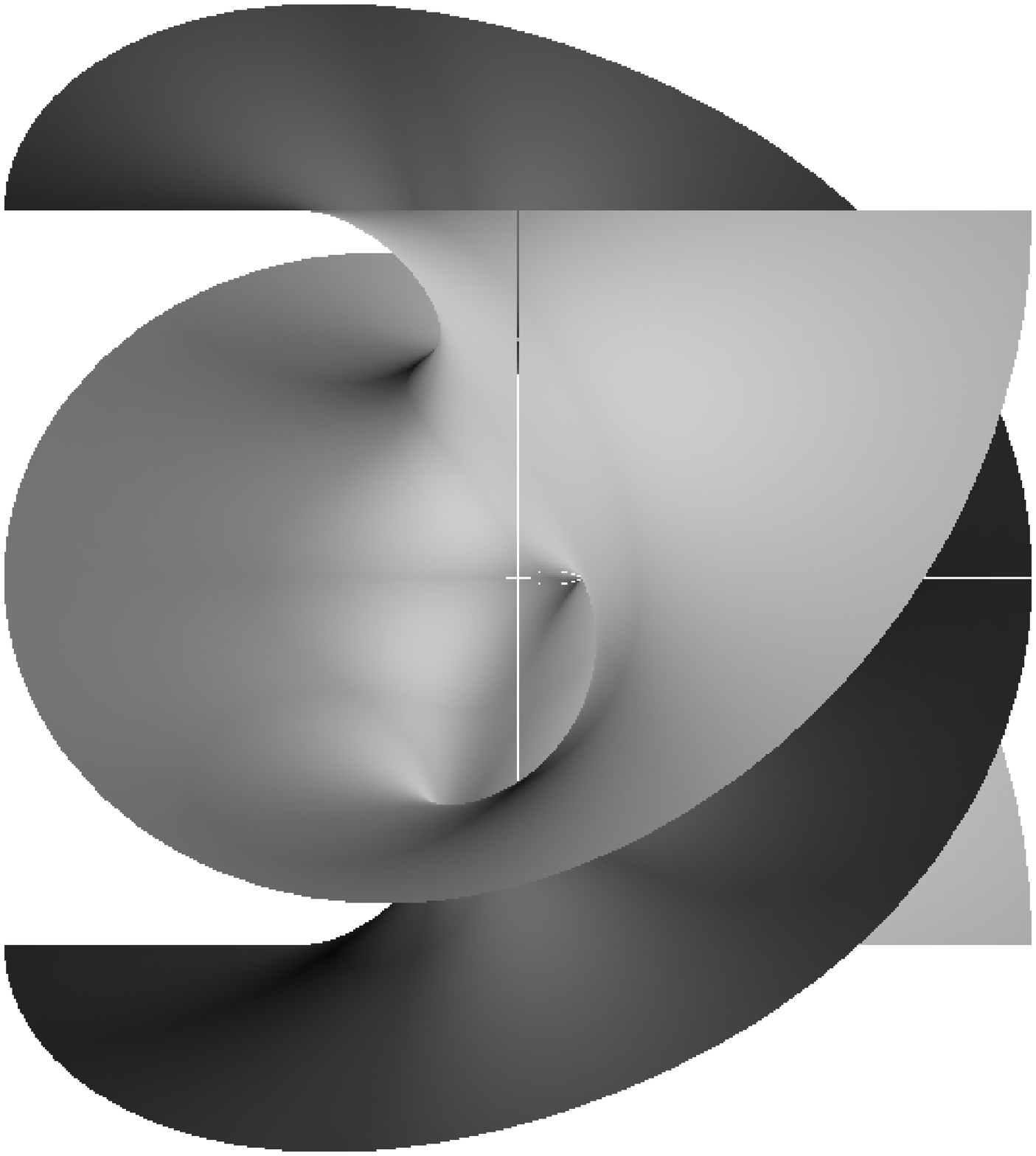} \hfill

\def\thecaption{
{\em Left Column}: The domain ${\cal D}$ with the image under $X$ of
$\partial {\cal D}$ above and  $\Omega$ the projection of $X({\cal D})$
below.  The curves of $X(\partial{\cal D})$ are labelled by their
preimages in $\partial {\cal D}$.  {\em Right Column}: The graph
$X({\cal D})$ at the top, with $X({\cal D})\cup X(\hat{\cal D})$ in
the middle.
 Extending the surface by rotation about the vertical axis
produces the fundamental domain, bottom right. Note that these surfaces
are tilted forward, so $x_3(C)$ does not appear to be monotonic, which it
actually is.}
\caption{\thecaption}
\endfig

We choose to integrate  from $O\in T^2$ in the Weierstrass 
representation (1.2) and we place $X(O)$ at
the origin of ${\Bbb R}^3$; i.e. $X (O) = {\bf O} := (0,0,0)$.
 With this normalization, it follows from Theorem 1 that 
the v-diagonal is mapped into the $x_3$-axis. 
The form of $dh$ in (1.7) and the fact that $z$ takes 
values of the form $it$, 
$t \leq 0$, on this axis, imply that
{\em $x_3$ increases monotonically as one descends the 
v-diagonal.}  

  The segment of the h-diagonal passing through $O$ 
and terminating at the punctures must be mapped onto a horizontal line through
$X(O) = {\bf O}$ in ${\Bbb R}^3$. Observing that $g(O) = g(O^\prime) = -1$,
 while $g = +1$ at the other fixed points of the normal
symmetry (where $z = \infty$), we can conclude that this horizontal line must
be the $x_2$-axis.  Orientation considerations or a direct computation of
$x_2$ imply that {\em $x_2$ is decreasing as one travels from $O$ along $H_1$ 
on the h-diagonal.}

\begin{proposition}\begin{enumerate} 

\item[(i)] $X (\interior{\cal D})$
is the bounded graph  of a function $F$ over an unbounded domain
$\Omega$ in the half plane $\{ (x_1, x_2)| x_1 \leq 0 \}$. The boundary
of $\Omega$ consists of the $x_2$-axis, together with a compact convex
curve $c$, that begins and ends at $(0, 0)$ and is symmetric with
respect to the $x_1$-axis. The function
$F$ satisfies  $-T/2 \leq F(x_1, x_2) \leq a$,
 where $T$ is the translational period
defined in (1.10), and $a > 0$ is defined by
$X(b) = (0, 0, a)$, where $b$ is the endpoint of $E$ in Figure~8;

\item[(ii)] $X (\partial {\cal D})$ is an embedding: $X(H_1)$ is the
negative $x_2$-axis; $X (H_2)$ is the positive $x_2$-axis translated 
 by $(0,0,-T/2)$; $X (E)$ is
the segment $(0, a)$ on the $x_3$-axis; $X (\hat E)$ is the segment
$(-a, -T/2)$ on the $x_3$-axis; $X (C)$ is a monotone graph over the curve
$c$ defined in (i) above. It
joins $X(b) = (0, 0, a)$ to $X(\hat b) =(0, 0, -a)$.
\end{enumerate}
\end{proposition}

We will prove Proposition 1 in the next section. The reader may
wish to read Section~3.3 first, where Theorem 3 is proved using
 Proposition 1.

\bigskip

\subsection{The proof of Proposition 1}

We will follow the boundary of ${\cal D}$, beginning at $O$ and we will
show 
 that $\partial {\cal D}$ is embedded by $X$ and
that---except for the vertical line segments  $X (E)$ and $X (\hat E)$ 
on the $x_3$-axis---$X (\partial {\cal D})$ is a graph over a curve in
the $(x_1, x_2)$-plane.  Along the way, we show that points of $ {\cal D}$
near $\partial {\cal D}$ are mapped by $X$ to points that 
 project into $\Omega$, and that points of ${\interior {\cal D}}$ near 
  $\partial {\cal D}$ are mapped by $X$ 
 to a graph over a neighborhood of $\partial \Omega$ in
$\interior{\Omega}$.  This, together with the behavior of $X$ near the
puncture, will allow us to conclude that $X (\interior{\cal D})$ is
a graph over $\interior{\Omega}$.

We define $$p : = \pi \circ X,$$ where
$\pi$ is the projection of ${\Bbb R}^3$ onto the $(x_1, x_2)$-plane.
Because $\vert g \vert >1$ on $\interior{\cal D}$, $p$ is an
immersion.  Beginning at $O \in T^2$,  which is mapped to ${\bf O} \in {\Bbb
R}^3$, the v-diagonal is mapped into the $x_3$-axis and $x_3$
increases monotonically as one descends the v-diagonal.\footnote{See
the computation of $dh(\dot E)$ near the end of Section 1.5.}  In 
particular,  $p(E) = (0, 0)$ and
$$X (b)= (0, 0, a)$$
for some $a >0$. We may describe the normal symmetry of the
surface by
$X\circ r_P =s\circ X,$ where
\begin{eqnarray}
s(x_1, x_2, x_3) = (x_1, -x_2, -x_3).
\end{eqnarray}
 Writing $\hat b := r\!_{_P}b$,
 we have
$$X (\hat b) = s\circ X(b) = (0, 0, -a).$$  We have shown
earlier that $X$ maps $H_1$, monotonically
onto the nonpositive $x_2$-axis. Therefore $p =X$ on $H_1$. 

 At $O \in T^2$, $g
=-1$.  Moreover,  $g$ is real and initially decreasing along $H_1$; after 
passing through
the vertical point, where $g =\infty$, $g$ continues to decrease until one
reaches a
branch point, where  $g$ has a value between $1$ and $\infty$. (We know
from Theorem~1 that $\lambda<1$, so the branch point of $g$ comes
before one reaches the puncture.) Then $g$
increases to $\infty$ as $H_1$ 
diverges to the end.  We can deduce two things from this.

First, because ${\cal D}$ is to the left of $H_1$ as one travels along
$H_1$ away from $O$, a neighborhood of $H_1$ in ${\cal D}$ is a
graph over a one-sided neighborhood of the negative $x_2$-axis, on the
side where $x_1 \leq 0$.  Consequently a suitably small neighborhood in
${\cal D}$ of the end, which is asymptotic to a helicoidal graph, must be a
graph over a region in $\{ (x_1, x_2) \vert x_1 \leq 0 \}$. 

Second, we can conclude that the surface normal rotates 
counter-clockwise as one
moves away from $O$ and keeps rotating in this sense until one
reaches the branch point.  The two straight lines on the surface, which
cross at $\bf O$,  are geodesic asymptotic curves. On the
 conjugate minimal immersion $X^\ast$, defined in a
neighborhood of $O \in T^2$, the lines correspond to geodesic principal
curvature lines, with principal curvatures given by the rate of change of
the surface normal along $E$ and $H_1$.  Since $O \in T^2$ is not a
branch point of $g$ and $X$ and $X^\ast$ have the same principal curvatures,
it follows that the normal to the surface rotates clockwise as one
descends $E$ from $O$. Because ${\cal D}$ is on the right as one descends
$E$, it follows that, near $X (E)$, $X ({\cal D})$ is a graph over a
region in $\{ (x_1, x_2) \vert x_1 < 0 \}$.  (This can also be deduced
from computation of the values of $g$ along $E$.)  It also implies that,
near  $X (b) = (0, 0, a)$,
$X (C - \{b\})$  is a graph over
a curve in $\{ (x_1, x_2) \vert x_1 < 0 \}$. 

Along $C$, $z$ is unitary.  On $C$, $z$ is branched only at the midpoint
of $C$, where $z = -e^{-i\rho}$.  We may use the phase of $z$ to
parameterize $C$ from $b$, where $z = -i$, to the midpoint of $C$, where
$z = -i e^{-i (\pi/2 + \rho)}$. That is 
$$C (t): = b + f (t) e^{i\pi/4}$$
is the point where $z \circ C (t) = -ie^{-it}$,  $0 \leq t \leq \pi/2 + \rho
< \pi$.

From Remark~2 at the end of Section~1, $x_3 \circ II(t)$ is a strictly
 decreasing function. (See Figure~7.) But $II$ is one half of $C$, from $b$
to the midpoint. 
By using the
normal symmetry
 we can conclude that $x_3$ is a decreasing function on all
of $C$.
In particular $X$  embeds $C$.

Because $g$ is unitary on $C$, the projection 
$$c := p (C) = \pi \circ X(C)$$ 
is a plane curve whose normal vector at $c(t)$ is $g(C(t))$.
The Gauss map is branched at $b$ and $\hat b$, where $g$ assumes complex
conjugate values, but has no branch points on the interior of $C$.  Hence
$c$, which begins and ends at $(0,0)$ has a normal that turns at a nonzero
rate all along $c$, and turns more than $\pi$ but less than $2\pi$ from
beginning to end.\footnote{$180^\circ$ rotation around the
v-diagonal produces a curve $\tilde C$ on $T^2$ with $g (\tilde C)
=-g (C)$.  Since the degree of $g$ is two, $g (C)$ is an arc of length
less than $2\pi$.  Because it begins and ends at $(0, 0)$ and has $g = +
1$ at its midpoint, $g (C)$ is an arc of length greater than $\pi$.
Alternatively, using the fact (established in Theorem~1) that $\rho>0$, we can 
  show that the $g(C)$ is an arc of length less than $2\pi$.} 
 We
 have already established that near $X (b)$, $X (C)$ projects into the
half space $\{ (x_1, x_2) \vert x_1 < 0 \}$. (See Figure~9.) Hence $c$ is an
embedded convex curve
in $\{ (x_1, x_2) \vert x_1 \leq 0 \}$ that begins and ends at $(0,0)$.  It is
symmetric with respect to the $x_1$-axis: Since $\pi \circ X = \hat s
\circ \pi$, where $ \hat s(x_1, x_2) := (x_1, -x_2)$, 
$$c(t) = p\circ C(t) = \pi \circ X(C(t)) = \pi\circ s\circ X(C(-t)) = 
{\hat s}\circ \pi \circ X(C(-t)) = {\hat s} (c(-t)).$$
We define $\Omega$ to be the unbounded region of $\{ (x_1, x_2) \vert x_1 \leq
0 \}$ bounded by the $x_2$-axis and $c$.  Observe that
 the normal to the surface along
$c$ points {\em out} from $\Omega$.  Since $c$ is strictly convex and the
Gaussian curvature $K$ is never zero 
along $X(c)$, except at the end points, the projection of
 any curve on the surface
orthogonal to $X(C)$ must lie in $\Omega$ for points sufficiently close to 
$X(C)$.  In particular, $\interior{C}$ has a neighborhood ${\cal U}$ with $p
({\cal U}) \subset \Omega$. Because ${\cal D}$ is to the right of $C$ as one
travels from $b$ toward $\hat b$, $X ({\cal U}\cap {\cal D})$ is {\em below}
 $X
(\interior{C})$,
 and $X ({\cal U}\cap \cal D)$ is a graph over $p ({\cal U} \cap
{\cal D}) \subset D$.

We can now continue around $\partial {\cal D}$.  $X (\hat E)$ is a segment of
the $x_3$-axis beginning at $X(b_2) = (0, 0, -a)$ and descending to
$X (O^\prime)$. Since  successive rotation about the horizontal lines
generates a translational period, $X( O') = (0, 0, -T/2)$, where $(0,
0, T)$ is the period vector of the singly periodic surface. 

 Near $X(\hat b)$, $X ({\cal D})$ projects into $\Omega$. 
 The normal to the surface along
$\hat V$ turns clockwise as one descends from $\hat b$ to $O^\prime$ and 
$g
(O^\prime) = -1$.  Since ${\cal D}$ lies to the right of $\hat E$ as one
ascends, points of ${\cal D}$ near $\hat E$ are mapped
 by $p = \pi \circ X$
to points near $(0, 0)$ lying in $\{ (x_1, x_2) \vert x_1 <
0\}$ and to the right of the tangent line to $c$ at its end point (i.e. to
the right of $g (\hat b)^\bot$).  In particular, near $\hat E$, ${\cal D}$ is a
graph over region of $\Omega$.

Arguments similar to those given above easily show  that $X (H_2)$ is the
ray $$\{(O, \tau, -T/2) | \tau >0 \},$$ and $\tau \to \infty$ as one approaches
the puncture at the end of $H_2$.  We already know that, near the puncture, 
$X ({\cal D})$ is a graph over a region of ${ \Omega}$.  We also know this
at $O^\prime$, the other end of $H_2$.  Since $g$ is never unitary on $H_2
- \{ O^\prime \}$, it follows that $H_2$ has a neighborhood in ${\cal D}$
that is a graph over a neighborhood in $\Omega$ of the positive $x_2$-axis.

We now have established statement (ii) as well as the fact that
 $X(\partial{\cal D})$ projects onto $\partial{\Omega}$. In
addition, we have shown that a neighborhood of
$\partial {\cal D} $ in ${\interior{\cal D}}$,  say ${\cal R}$, 
is mapped to a
 graph over a neighborhood in $\Omega$ of $ \partial \Omega$ and that a
neighborhood in ${\cal D}$ of the 
puncture is a graph over a region in $\Omega$ of the form $\{ \vert(x_1, x_2)
\vert > M\}$, for $M$ sufficiently large.  This means that $X({\cal R}
-\partial {\cal D})$ is a graph over a region in $\Omega$ of the form 
$p({\cal R})\cup \{ (x_1,x_2) \in \Omega \vert ~\vert (x_1, x_2) \vert > M \}$. 
 The complement of this
region in $\Omega$ is compact and simply connected and  has a boundary over 
which $X ({\partial {\cal R}})$
is graph.  Since $p$ is an immersion on
${\cal D}$, it follows that $X({\cal D} - {\cal R})$ is a graph over this
complementary region. Hence  $X (\interior{\cal D})$ is a graph of some 
function, which we call $F$,  over
$\interior{\Omega}$, and is asymptotic to a helicoidal graph as $\vert (x_1,
x_2) \vert \to \infty$. In particular, we have proved the
first claim of statement (i). Because $\Omega$ lies in the half plane $\{ (x_1,
x_2) \vert x_1 \leq 0 \}$ and $-T/2 \leq x_3 \leq a$ on $X (\partial
{\cal D})$, it follows that $X ({\cal D})$ lies in the slab of ${\Bbb
R}^3$ defined by these two constraints. This gives the bounds on $F$
stated in statement (i).

\subsection {The proof of Theorem 3}

It suffices to prove that a fundamental piece of the surface, modulo
translations, is embedded.
 We know that such a fundamental piece
 is made up of four copies of the closure of the graph $X ({\cal D})$,
 described in Proposition 1.
 Let $\hat {\cal D} = {{r\!_{_P}}} {\cal
D}$.

\begin{claim}
$X$ is an embedding of ${\cal D} \cup \hat {\cal
D}$ with values in the slab $\{ (x_1, x_2, x_3) \vert\; x_1 \leq 0, -T/2\leq
x_3 \leq T/2\}$ and boundary values consisting of vertical line segments 
and rays on $\{ x_1 = 0\}$. 
\end{claim}

Assuming the claim, which we will prove below, we  can 
complete the proof of  Theorem 3.

 The rest of
a fundamental domain of the surface is produced from $X ({\cal D} \cup
\hat {\cal D})$ by $180^\circ$ rotation about the $x_3$-axis; according
to (1.14), $X\circ\mu_{vert}(p) = (-x_1,-x_2,x_3 )(p)
=:\sigma\circ X(p)$

This means that $\sigma \circ X
({\cal D} \cup \hat {\cal D})$ lies in the slab
 $\{ (x_1, x_2, x_3) \vert\,
x_1 \geq 0$, $\vert x_3 \vert \leq T/2 \}$, and is therefore disjoint from $X
({\cal D} \cup \hat {\cal D})$ except along their common boundary, 
where we already know that $X$ is one-to-one.  Hence
$X$ is an embedding.

\begin{proof}[Proof of the Claim]
Recall from (3.1) that
$s (x_1, x_2, x_3) =
(\hat s (x_1, x_2), -x_3)$, where $\hat s(x_1, x_2):= (x_1, -x_2)$.   
Since $X (\hat {\cal D}) =X (r_P \circ {\cal D})=
 s \circ X ({\cal
D})$, $X (\interior{\hat {\cal D}})$ is the graph of the function
\begin{equation}
\hat F (x_1, x_2) := -F \circ \hat s (x_1, x_2)
\end{equation}
on $\Omega -\{ {\bf O}\}$, where $F$ is the function, 
described in Proposition 1,  whose graph
 is $X (\interior{\cal D})$.  Since $\hat s$ interchanges the positive and
negative $x_2$-axis, on which $F = -T/2>0$ and $F = 0$, respectively, (3.2)
implies that $\hat F > F$ on these rays. 

 For $\vert (x_1, x_2)\vert$
large,  $X ({\cal D})$ is asymptotic to the graph of a helicoid over $\{
(x_1, x_2) \vert x_1 \leq 0, (x_1, x_2) \ne (0, 0) \}$ whose boundary
values are $-T/2$ on \{$x_2 > 0$\} and $0$ on \{$x_2 < 0$\}. 
From (3.2) it follows that $\hat F > F$ for $\vert (x_1,
x_2)\vert$ large.

Recall that $c$ is the projection of $X (C) $ onto 
the $(x_1, x_2)$-plane. 
Because $s\circ X (C) =X (C)$ and $c \circ \hat s  = c$, $\hat F = F$ on $c$.
According to Proposition 1, $X$ is on-to-one on $C$. Therefore $X$ is one-to-one
on a neighborhood of $C$, so  $F \neq \hat F$  near,
 but not on, $c$. 
 The line segments $X (H_1)$ and $X (H_2)$ are disjoint segments of
the $x_3$-axis and $s \circ X (H_1)$ and $s \circ X (H_2)$ meet $X (H_1)
\cup X (H_2)$ only at $X (b)$ and $X (\hat b)$, near which $X$ is an
embedding.

 We have shown at each point $(x_1, x_2)\in \partial {\Omega}$
 that either $\hat F(x_1, x_2) > F(x_1, x_2)$ or
$\hat F(x_1, x_2) = F(x_1, x_2)$ and that 
$\hat F >F$ for nearby points in $\interior{\cal D}$.
 Also, $\hat F(x_1, x_2) > F(x_1, x_2)$ 
for
 $|(x_1, x_2)|$   sufficiently large.  By the Maximum
Principle, $\hat F > F$ on $\interior{\Omega}$ because the set of
points in $\interior{\Omega}$ where $\hat F\leq F$ is bounded and can have
no limit points on $\partial \Omega$.  
Thus the graphs of $\hat F$ and $ F$ over $\interior{\Omega}$ are
disjoint, which means 
that $X ({\cal D} \cup \hat {\cal D})$
is an embedding, as claimed. The boundary of $X ({\cal D} \cup \hat {\cal D})$
 consists of the two
vertical line segments of $X ({\cal D})$ and $X (\hat {\cal
D})$, together with the two horizontal rays $\{ x_3 = -T/2,\; {x_1 =
0,\; x_2 \geq 0} \}, \{ x_3 = T/2,\;  {x_1 = 0,\; x_2 \leq
0}\}$ and the  horizontal line $\{ x_1 = x_3 = 0 \}$.
\end{proof}

\bibliographystyle{plain}

\end{document}